\documentclass{aptpub}
\usepackage{graphicx}
\usepackage{amsmath}
\usepackage{amsfonts}
\usepackage{amssymb}
\usepackage{setspace}
\usepackage{cite}

\renewcommand{\crnotice}{}

\renewcommand{\AAP}{\emph{Adv. Appl. Probab.\ }}

\newcommand{\be}{\begin{equation}}
\newcommand{\ee}{\end{equation}}
\newcommand{\bea}{\begin{eqnarray}}
\newcommand{\beas}{\begin{eqnarray*}}
\newcommand{\no}{\nonumber}
\newcommand{\eea}{\end{eqnarray}}
\newcommand{\eeas}{\end{eqnarray*}}

\newcommand{\selfnote}[1]{\textbf{
    \\ ***Note to Self:} {\textsf{#1}} \textbf{***} \\}
\newcommand{\qed}{\hfill $\Box$}

\newcommand{\remove}[1]{}

\newcommand{\resp}[1]{\,(\text{resp.}\,#1)}

\newcounter{cnt1}

\newcounter{cnt3}
\newcommand{\blr}{\begin{list}{$($\roman{cnt1}$)$}
 {\usecounter{cnt1} \setlength{\topsep}{0pt}
 \setlength{\itemsep}{0pt}}}
\newcommand{\bla}{\begin{list}{$($\betaph{cnt2}$)$}
 {\usecounter{cnt2} \setlength{\topsep}{0pt}
 \setlength{\itemsep}{0pt}}}
\newcommand{\bln}{\begin{list}{$($\arabic{cnt3}$)$}
 {\usecounter{cnt3} \setlength{\topsep}{0pt}
 \setlength{\itemsep}{0pt}}}
\newcommand{\el}{\end{list}}

\def\rar{\rightarrow}

\newcommand{\lam}{\lambda}
\newcommand{\Lam}{\Lambda}

\newcommand{\sg}{\sigma}
\newcommand{\1}{{\bf 1}}

\def\Ph{\Phi}

\def\mR{\mathbb{R}}
\def\mZ{\mathbb{Z}}
\def\mN{\mathbb{N}}

\def\mE{\mathbb{E}}
\def\mM{\mathbb{M}}

\authornames{B. Blaszczyszyn , D. Yogeshwaran}
\shorttitle{Ordering of random measures and shot-noise fields}

\begin{document}


\title{Directionally Convex Ordering of Random Measures,
Shot Noise Fields  and Some Applications to
Wireless Communications}\\[2ex]
\authorone[INRIA/ENS and Math. Inst. University of
  Wroc{\l}aw]{Bart\L omiej  B\L aszczyszyn\vspace{-5ex}}
\addressone{ENS DI TREC, 45 rue d'ulm, 75230 Paris, FRANCE.}
\noindent\authortwo[INRIA/ENS]{D. Yogeshwaran}
\addresstwo{ENS DI TREC, 45 rue d'ulm, 75230 Paris, FRANCE.}
\vspace{-4ex}

%


\begin{abstract}
Directionally convex ($dcx$) ordering  is a
tool for comparison of dependence structure
of random vectors that also takes into account the variability of the
marginal distributions.
When extended to random fields it concerns
comparison of all finite dimensional distributions.
Viewing locally finite measures as non-negative fields of
measure-values indexed by the bounded Borel subsets of the space,
in this paper we  formulate and study the
$dcx$ ordering of random measures on locally compact spaces.
We show that the $dcx$ order is preserved under some
of the natural operations considered on random measures and point
processes, such as deterministic displacement of points,
independent superposition and thinning as well as independent, identically
distributed marking. Further operations such as
position dependent marking and displacement of points
are shown to preserve the order on Cox point processes. We also
examine the impact of $dcx$ order on the second moment properties, in
particular on clustering  and on Palm distributions.
Comparisons of Ripley's functions, pair correlation functions
as well as   examples seem to
indicate that point processes higher in $dcx$ order cluster more.
\\
As the main result, we show that non-negative
integral shot-noise fields with respect to
$dcx$ ordered random measures inherit this ordering from the measures.
Numerous  applications of this result are shown,
in particular to comparison of various Cox processes and some
performance measures of wireless networks, in both of which shot-noise
fields appear as key ingredients.
We also mention a few pertinent open questions.
\end{abstract}


\section{Introduction}
\label{sec:Intro}

\remove{
\selfnote{This is a temporary introduction}

Thoughts
\begin{itemize}
\item Various generalizations of convexity in multi-dimension,
\item Ross conjecture,
\item Strong order of point processes,
\item Why Cox, why SN?
\item $dcx$ and clustering.
\item Completion of Miyoshi.
\end{itemize}
}

Point processes (p.p.) have been at the centre of various studies in
stochastic geometry, both theoretical and applied.
Most of the work involving quantitative analysis of p.p. have dealt with
Poisson p.p.. One of the main reasons being that characteristics of
Poisson p.p. are amenable to computations and yield nice closed form
expressions in many cases. Computations have been difficult
in great many cases, even for  Cox (doubly stochastic Poisson) p.p..

\paragraph{Comparison of point processes}
To improve upon this situation, qualitative,
comparative studies of p.p. have emerged as useful tools.
The first method of comparison of p.p. has been coupling or stochastic
domination (see \cite{Kamae77,Lindvall92,Rolski91}). In our
terminology, these are known as strong ordering of
p.p.. When two p.p. can be coupled, one turns out to be a subset of the other.
This ordering is very  useful for obtaining various bounds and proving
limit theorems.  However, using it one cannot
compare two different p.p. with same mean measures. An
obvious example is an  homogeneous Poisson p.p. and a stationary Cox
p.p. with the same intensity.
The question arises of what ordering is
suitable for such p.p.?
This is an important question since
it is expected that by comparing p.p. of the  same intensity one should
achieve a tighter bound than by coupling.
For some more details on strong ordering of p.p. and need for other orders,
see remarks in~\cite[Section 5.4 and Section 7.4.2]{Muller02}.

\paragraph{From convex to {\em  dcx} order}
Two random variables $X$ and $Y$ with the same mean $\sE(X)=\sE(Y)$
can  be compared by how "spread out" their distributions are.
This {\em statistical variability} (in statistical ensemble)
is captured to a limited extent by the variance, but more fully
by convex ordering, under which $X$ is less than $Y$ if and only if
for all convex $f$, $\sE(f(A))\le\sE(f(B))$.
In multi-dimensions, besides different statistical
variability of marginal distributions, two random vectors can
exhibit different {\em dependence properties} of their coordinates.
The most evident example here is comparison of the  vector
composed of several copies of one random variable to a vector
composed of independent copies sampled from the same distribution.
A useful tool for comparison of the dependence structure of random
vectors with fixed marginals is the {\em supermodular order}.
The $dcx$ order is another integral order (generated by a class of
$dcx$ functions in the same manner as convex functions generate the
convex order) that can be seen as  a generalization of the supermodular one,
which in addition takes into account the variability of the marginals
(cf~\cite[Section 3.12]{Muller02}).
It can be naturally extended to random fields by
comparison of all finite dimensional distributions.

\paragraph{The {\em dcx} order of random measures}
In this paper we make an obvious further extension that
consists in $dcx$ ordering of locally finite  measures
(to which belong p.p.) viewed as  non-negative fields of
measure-values on all bounded subsets of the space. We show that the
$dcx$
order is preserved  under some of the natural operations considered on
random measures and point processes, such as independent superposition
and thinning. Also, we examine
the impact of $dcx$ order on the second moment properties, in
particular on clustering,  and Palm distributions.

\paragraph{Integral shot-noise fields}
Many interesting characteristics of random measures, both in the
theory and in applications  have the form of integrals of some
non-negative kernels. We call them {\em integral  shot-noise fields}.
For example, many classes of Cox p.p., with the most general being
L\'{e}vy based Cox p.p. (cf.~\cite{Hellmund08}),
have stochastic intensity fields, which are
shot-noise fields. They are also key ingredients of
the recently proposed, so-called ``physical''  models
 for  wireless networks, as we will explain in
what follows (see also~\cite{GK00, Baccelli01,Dousse_etal2006}).
It is thus particularly appealing to study the shot-noise fields
generated by $dcx$ ordered random measures.

Since integrals are linear operators on the space of measures,
and knowing that a linear function of a vector is trivially  $dcx$,
it is naturally to expect that the integral shot-noise fields with respect to
$dcx$ ordered random measures will inherit this ordering from the measures.
However, this property cannot be concluded immediately from the finite
dimensional $dcx$  ordering of measures.
The formal proof of this fact that is the main result of this paper
involves some arguments from the
theory of integration combined with the closure property of $dcx$
order under joint weak convergence and convergence in~mean.

\paragraph{Ordering in  queueing theory and  wireless communications}
The theory of stochastic ordering provides
elegant and efficient tools for comparison of random objects
and is now being used in many fields.
In particular in queueing theory context, in~\cite{Ross78},
Ross made a conjecture
that replacing a stationary Poisson arrival process in a single server
queue by a stationary  Cox p.p. with the same intensity  should increase the
average customer delay.
There have been many variations
of these conjectures which are now known as Ross-type
conjectures. They  triggered the interest in comparison of
queues with similar inputs (\cite{Chang91,Miyoshi04a,Rolski89}).
The notion of a $dcx$ function was partially developed and used
in conjunction  with the proving of Ross-type
conjectures (\cite{LMeester93,LMeester99,Shaked90}). Much earlier to
these works, a comparative study of queues motivated by neuron-firing models can
be found in~\cite{Huffer84a}. Also comparison of variances of point processes and
fibre processes was studied in \cite{Stoyan83} and hence it can
be considered as a forerunner to our article.
The applicability of these results has
generated sufficient interest in the theory of stochastic ordering
as can be seen from the diverse results in the book of M\"uller and
Stoyan (\cite{Muller02}).
As most works on ordering of
p.p. were motivated by applications to queueing theory, results were
primarily focused on one-dimensional point processes.
An attempt to rectify the lack of
work in higher dimensions was made in \cite{Miyoshi04}, where
comparison results for shot-noise fields of spatial stationary Cox
p.p. were given. The results of~\cite{Miyoshi04} are the starting
point of our investigation.

Our interest in ordering of point processes, and in particular in the
shot-noise fields they generate,
has roots in the analysis of  wireless communications, where these objects
are primarily used to model the so called {\em interference}
that is the total power
received from many emitters scattered in the plane or space
and sharing the common Hertzian medium.
According to a new emerging methodology,
the interference-aware stochastic
geometry modeling of wireless communications provides a way of
defining and computing macroscopic properties of large wireless
networks  by some averaging over all
potential random patterns for node locations in
an infinite plane and   radio channel
characteristics, in the same way as queuing theory provides
averaged response times or congestion over all potential arrival
patterns within a given parametric class.
These macroscopic
properties will allow one to characterize the key dependencies of the
network performance characteristics
in function of a relatively small number of parameters.

In the above context, Poisson distribution of emitters/receiver/users
is  often too simplistic. Statistics show that the real patterns
of users exhibits more clustering effects (``hots spots'') than observed in
an homogeneous Poisson point processes. On the other hand,  good
packet-collision-avoidance mechanisms  scheme should create
some ``repulsion'' in the pattern of nodes allowed to access
simultaneously to the channel. This rises questions about the analysis
of non-Poisson models, which could be to some extent tackled  on the
ground of the theory of stochastic ordering.
Interestingly, we shall  show that there are certain performance
characteristics  in wireless networks that
{\em improve} with more variability in the input process.

\paragraph{The remaining part of the article is organized as follows.}
In the next section,
we will
present the main  definitions and state the main results concerning
$dcx$ ordering of the integral shot-noise fields.
Section~\ref{sec:ord_rand_meas} will explore the various
consequences of ordering of random measures. The proofs of the main
results are given in
Section~\ref{sec:OSN}. Examples illustrating the use and application of the theorems shall be presented in Section~\ref{sec:examples}.
Section~\ref{sec:appln} will sketch some of the possible applications of results in the context of wireless communications.
Finally, we conclude with some remarks and questions in
Section~\ref{sec:concl}. There is an Appendix (Section~\ref{sec:Appendix})
containing some properties of stochastic orders and their extensions
that are used in the paper.
\section{Definitions and the Main Result}
\label{sec:defn_main}

The order $\leq$ on $\mR^n$ shall denote the component-wise partial
order, i.e., $(x_1,\ldots,x_n) \leq (y_1,\ldots,y_n)$ if $x_i \leq
y_i$ for every $i$.

\begin{defn}
\label{defn:dcx_fn}
\begin{itemize}
\item We say that a function $f:\mR^d \rar R$ is
{\em directionally convex}~($dcx$) if for every $x,y,p,q \in \mR^d$ such
that $p \leq x,y \leq q$ and $x+y = p+q$,
\[f(x) + f(y) \leq f(p) + f(q). \]
\item Function $f$ is said to be {\em directionally concave}~($dcv$)
  if the inequality in
the above equation is reversed.
\item  Function $f$ is
said {\em directionally linear} ($dl$) if it is $dcx$ and $dcv$.
\end{itemize}
\end{defn}

Function $f=(f_1,\ldots,f_n):\mR^d
\rar \mR^n$ is said to be
$dcx$($dcv$) if each of its component $f_i$ is $dcx$($dcv$).
Also, we shall abbreviate {\em increasing} and $dcx$ by $idcx$ and
{\em decreasing}
and $dcx$ by $ddcx$. Similar abbreviations shall be used for $dcv$
functions.
Moreover, we abbreviate {\em non-negative} and $idcx$ by  $idcx^+$.

In the following, let $\mathfrak{F}$ denote some class of
functions from $\mR^d$ to  $\mR$. The dimension $d$ is assumed to be
clear from the context. Unless mentioned, when we state
$\sE(f(X))$ for $f \in \mathfrak{F}$ and $X$ a random vector, we
assume that the expectation exists, i.e., for each random vector $X$
we consider the sub-class of $\mathfrak{F}$ for which the
expectations exist with respect to (w.r.t) $X$.

\begin{defn}
\label{defn:dcx_order}
\begin{itemize}
\item Suppose $X$ and $Y$ are real-valued random vectors of the
same dimension. Then  $X$ is said to be {\em less than $Y$ in
$\mathfrak{F}$ order} if
$\sE(f(X)) \leq \sE(f(Y))$ for all $f \in \mathfrak{F}$
(for which both expectations are finite). We shall
denote it as $X \leq_{\mathfrak{F}} Y$.

\item Suppose $\{X(s)\}_{s \in S}$ and $\{Y(s)\}_{s \in S}$ are
real-valued random fields, where  $S$ is an arbitrary index set.
We say that $\{X(s)\} \leq_{\mathfrak{F}}
\{Y(s)\}$ if for every $n \geq 1$ and $s_1,\ldots,s_n\in S$,
$(X(s_1),\ldots,X(s_n)) \leq_{\mathfrak{F}}(Y(s_1),\ldots,Y(s_n)).$

\end{itemize}
\end{defn}
In  the remaining part of the paper, we will mainly consider
$\mathfrak{F}$  to be the class of $dcx$, $idcx$ and $idcv$ functions;
the negation of these functions give rise
to $dcv,ddcv$ and $ddcx$ orders respectively.
If $\mathfrak{F}$ is the class of increasing
functions, we shall replace $\mathfrak{F}$ by $st$ (strong) in the
above definitions. These are standard notations used in literature.

As concerns random measures, we shall work in the set-up of \cite{Kallenberg83}.
Let $\mE$ be a
locally compact, second countable Hausdorff (LCSC) space. Such spaces are
polish, i.e., complete and separable metric space. Let
$\texttt{B}(\mE)$ be the Borel $\sg$-algebra and $\texttt{B}_b(\mE)$
be the $\sg$-ring of {\em bounded, Borel subsets} (bBs). Let $\mM =
\mathbb{M(\mE)}$
be the space of non-negative Radon measures on $\mE$. The Borel
$\sg$-algebra $\mathcal{M}$ is generated by the mappings $\mu
\mapsto \mu (B)$ for all $B$ bBs. A random measure
$\Lambda$ is a mapping from a probability space $(\Omega,\mathcal{F},\sP)$ to
$(\mM,\mathcal{M})$. We shall call a random measure $\Phi$ a p.p. if
$\Phi \in \bar{\mN}$, the subset of counting measures in $\mM$. Further,
we shall say a p.p. $\Phi$ is simple if a.s. $\Phi(\{x\}) \leq 1$ for all $x \in \mE$.
Throughout, we shall use $\Lambda$ for an arbitrary random measure
and $\Phi$ for a p.p.. A random measure $\Lambda$ can be
viewed as a random field $\{\Lambda(B)\}_{B \in \texttt{B}_b(\mE)}.$
With this viewpoint and the previously introduced
notion of ordering of random fields, we define ordering of
random measures.

\begin{defn}
\label{defn:dcx_rm} Suppose $\Lambda_1(\cdot)$ and $\Lambda_2(\cdot)$ are
random measures on $\mE$. We say that $\Lambda_1(\cdot) \leq_{dcx}
\Lambda_2(\cdot)$ if for any $I_1, \ldots, I_n$ bBs in
$\mE$,
\begin{equation}
\label{condn:dcx_rm}
 (\Lambda_1(I_1),\ldots,\Lambda_1(I_n)) \leq_{dcx} (\Lambda_2(I_1),\ldots,\Lambda_2(I_n)).
\end{equation}
\end{defn}
The definition is similar for other orders, i.e., when $\mathfrak{F}$ is the class of $idcx/idcv/ddcx/ddcv/st$ functions.

\begin{defn}\label{defn:ISN}
Let $S$ be any set and $\mE$ a LCSC space. Given a random measure
$\Lambda$ on $\mE$ and a measurable (in the first variable alone)
{\em response function} $h(x,y) :\mE \times S \rar \bar{\mathbb{R}}^+$
where $\bar{\mathbb{R}}^+$ denotes the completion of positive
real-line  with infinity, the {\em (integral) shot-noise field} is
defined as
\begin{equation}
\label{eqn:sn_rm} V_{\Lambda}(y) = \int_{\mE}h(x,y)\Lambda(dx).
\end{equation}
\end{defn}

With this brief introduction, we are ready to state our key result
that will be proved in Section~\ref{sec:ISN}.
\begin{thm}
\label{thm:isn_rm}
\begin{enumerate}

\item \hspace{0.3em} If $\Lambda_1  \leq_{idcx\resp{idcv}} \Lambda_2$,
then $\{V_{\Lambda_1}(y)\}_{y \in S} \leq_{idcx\resp{idcv}}
\allowbreak
\{V_{\Lambda_2}(y)\}_{y \in S}$.

\item Let $\sE(V_{\Lambda_i}(y)) < \infty$,  for all $y \in S$,
  $i=1,2.$ If $\Lambda_1  \leq_{dcx} \Lambda_2$,
then $\{V_{\Lambda_1}(y)\}_{y \in S} \leq_{dcx} \{V_{\Lambda_2}(y)\}_{y \in S}$.

\end{enumerate}
\end{thm}


The first part of the above theorem
for the one-dimensional marginals of bounded shot-noise fields
generated by lower semi-continuous response functions
is proved in~\cite{Miyoshi04}
for the special case of spatial stationary Cox p.p..
It is conspicuous that we have
generalized the earlier result to a great extent.
This more general result will be used in many places in this paper,
in particular to prove ordering of  independently, identically  marked
p.p. (Proposition \ref{prop:mpp_idcx}),
Ripley's functions (Proposition~\ref{prop:Ripley}),
Palm measures (Proposition~\ref{prop:palm_pp_idcx}),
independently marked Cox processes (Proposition~\ref{prop:mark_cox}), extremal
shot-noise fields (Proposition~\ref{prop:lo_msn}).
Apart form these results, Sections~\ref{sec:examples}
and~\ref{sec:appln}
shall amply demonstrate examples and applications that shall need
Theorem~\ref{thm:isn_rm}.



\section{Ordering of Random Measures and Point Processes}
\label{sec:ord_rand_meas}
We shall now give a sufficient condition for random measures to be
ordered, namely that
the condition (\ref{condn:dcx_rm}) in Definition \ref{defn:dcx_rm}
needs to be verified only for disjoint bBs.
The necessity is trivial.
This is a
much easier condition and will be used many times in the remaining
part of the paper.

\begin{prop}
\label{prop:char_dcx_rm} Suppose $\Lambda_1(\cdot)$ and $\Lambda_2(\cdot)$
are two random measures on $\mE$. Then $\Lambda_1(\cdot)
\leq_{dcx}\Lambda_2(\cdot)$ if and only if condition
(\ref{condn:dcx_rm}) holds for all {\em mutually disjoint} bBs.
The same results holds true for $idcx$ and $idcv$ order.
\end{prop}

\begin{proof} We need to prove the 'if' part alone. We shall prove for $dcx$ order and the same argument is
valid for $f$ being $idcx$ or $idcv$. Let condition (\ref{condn:dcx_rm}) be satisfied for all mutually  disjoint bBs. Let
 $f:\mR^n_{+} \rar R$ be  $dcx$ function and $B_1,\ldots,B_n$ be bBs. We can choose mutually disjoint
bBs $A_1,\ldots,A_m$ such that $B_i = \cup_{j \in J_i} A_j$ for all $i$. Hence $\Lambda(B_i) = \sum_{j \in
J_i} \Lambda(A_j).$ Now define $g:\mR^m_+ \rar \mR^n_+$ as $g(x_1,\ldots,x_m) = (\sum_{j \in J_1}x_j,\ldots,\sum_{j \in
J_n}x_j).$ Then $g$ is $idl$ and  so $f \circ g$ is $dcx$. Moreover,
$f(\Lambda(B_1),\ldots,\allowbreak \Lambda(B_n)) =
f \circ g(\Lambda(A_1),\ldots,\Lambda(A_m))$ and thus the result for $dcx$ follows. \qed
\end{proof}

\remove{
Now we turn to case of $lo$ order.we use notation similar to above. Let $a_i \in \mR^+, 1 \leq i \leq n.$
\begin{eqnarray*}
\sP(\Lam_1(B_i) \leq a_i, \forall i) & = & \sP(\sum_{j \in J_i}\Lam_1(A_j) \leq a_i \forall i) \\
& = & \int_{\sum_{j \in J_1} x_{j} \leq a_1} \ldots \int_{\sum_{j \in J_n} x_{j} \leq a_n} \sP(\Lam_1(A_j) \leq x_j \forall j) \prod_{j} dx_j \\
& \geq &   \int_{\sum_{j \in J_1} x_{j} \leq a_1} \ldots \int_{\sum_{j \in J_n} x_{j} \leq a_n}\sP(\Lam_2(A_j) \leq x_j \forall j) \prod_{j} dx_j \\
& = &\sP(\Lam_2(B_i) \leq a_i, \forall i)
\end{eqnarray*}
}

\subsection{Simple Operations Preserving Order}
\label{sec:simple_opns}

Point processes are special cases of random measures and as such will
be subject to
the considered ordering. It is known  that each p.p. $\Phi$
on a LCSC space $\mE$ can be represented as a countable sum
$\Phi=\sum_{i}\varepsilon_{X_i}$ of Dirac measures
($\varepsilon_x(A)=1$ if $x\in A$ and 0 otherwise)
in such a way that $X_i$ are random elements in $\mE$.
We shall now show that all the three orders  $dcx,idcx,idcv$ preserve
some simple operations on random measures and p.p., as
deterministic mapping, independent identically distributed (i.i.d.)
thinning and independent  superposition.

Let $\phi:\mE\to\mE'$ be a measurable mapping to some LCSC space
  $\mE'$. By the image of a (random) measure $\Lambda$ by $\phi$ we
understand $\Lambda'(\cdot)=\Lambda(\phi^{-1}(\cdot))$.
Note that the image of a p.p. $\Phi$ by $\phi$ consists in
deterministic displacement of all its points by $\phi$.

Let  $\Phi=\sum_{i}\varepsilon_{x_i}$.
By {\em i.i.d. marking} of $\Phi$, with marks in some
LCSC space $\mE'$, we understand
a p.p. on the product space $\mE\times\mE'$, with the usual product
Borel $\sigma$-algebra, defined by
$\tilde{\Phi}=\sum_{i}\varepsilon_{(x_i,Z_i)}$, where
$\{Z_i\}$ are i.i.d. random variables (r.v.), so called marks,
on  $\mE'$.
By {\em i.i.d. thinning} of $\Phi$, we understand
$\overline{\Phi}=\sum_{i}Z_i\varepsilon_{x_i}$, where $Z_i$ are i.i.d.
0-1 Bernoulli random variables r.v.. The probability $\sP\{Z=1\}$ is called
the {\em retention probability}.
{\em Superposition} of p.p. is understood
as  addition of  (counting) measures. Measures on Cartesian products of LCSC
spaces are always considered with their corresponding product Borel $\sigma$-algebras.
\begin{prop}
\label{prop:mpp_idcx}
Suppose $\Lambda_i, i=1,2$ are random measures and $\Phi_i, i =1,2$ are p.p.. Assume that
$\Lambda_1\leq_{dcx\resp{idcx;\,idcv}}\Lambda_2$ and  $\Phi_1 \leq_{dcx\resp{idcx;\,idcv}} \Phi_2$.
\begin{enumerate}
\item Let $\Lambda'_i$  be the image of $\Lambda_i$, $i=1,2$,
by some mapping $\phi:\mE\to\mE'$.
Then
  $\Lambda'_1\leq_{dcx\resp{idcx;\,idcv}}\Lambda'_2$. As a special
  case, the same holds true for the displacement of points of
  $\Phi_i$'s by $\phi$.
\item Let $\Phi_i,i=1,2$, be simple p.p. and $\tilde{\Phi}_i, i =1,2$, be the corresponding
  i.i.d. marked p.p. with the same distribution of marks.
Then $\tilde{\Phi}_1 \leq_{dcx\resp{idcx\,;idcv}} \tilde{\Phi}_2$.
\item Then $\overline{\Phi}_i$ be i.i.d. thinning of $\Phi_i$, $i
  =1,2$,
with the same retention probability.
Then $\overline{\Phi}_1 \leq_{dcx\resp{idcx\,;idcv}} \overline{\Phi}_2$.

\item Let $\Lam_1'$ and $\Lam_2'$  be two random measures
such that $\Lam_1' \leq_{dcx\resp{idcx;\,idcv}} \Lam_2'$.  Assume that
$\Lam_i'$'s are independent of $\Lam_i$'s.
Then  $\Lam_1+\Lam_1' \leq_{dcx\resp{idcx;\,idcv}} \Lam_2 + \Lam_2'$,
where $+$ is understood as the addition of measures.
\item Suppose the random measures are on the product space  $\mE\times\mE'$.
Then
$\Lambda_1(\mE\times\cdot)\leq_{dcx\resp{idcx;\,idcv}}\Lambda_2(\mE\times\cdot)$,
provided the respective projections are Radon measures.

\end{enumerate}
\end{prop}
\begin{proof}
\begin{enumerate}
\item
The result follows immediately from the Definition~\ref{defn:dcx_rm}.

\item We shall prove $\tilde{\Phi}_1 \leq_{dcx} \tilde{\Phi}_2$ and the proof for the other
orders is similar. Since $\mE$ is a LCSC space, there exists a {\em null-array of partitions}
$\{B_{n,j} \subset \mE \}_{n \geq 1, j \geq 1}$, i.e., $\{B_{n,j}\}_{j
  \geq 1}$ form a finite partition of $\mE$ for every $n$
and $\max_{j \geq 1}\{|B_{n,j}|\} \to 0$ as $n \to \infty$ where
$|\cdot|$ denotes the diameter in any fixed metric
(see \cite[page~11]{Kallenberg83}).  For every $x \in \mE$, let
$j(n,x)$ be  the unique index such that $x \in B_{n,j(n,x)}$.
Let $\overline{Z} = \{Z_{n,j}\}_{n \geq 1, j \geq 1}$ be a family of
$\mE'$-valued i.i.d. random variables with distribution $F(\cdot)$.
Define marked p.p. $\tilde{\Phi}^n_i = \sum_{X_k \in \Phi_i}\varepsilon_{(X_k,Z_{n,j(n,X_k)})}$
for $i =1,2$. We shall now verify that the sequences $\tilde{\Phi}^n_i$'s satisfy the assumption of Lemma \ref{lem:weak_conv} with limits
$\tilde{\Phi_i}$'s respectively.

Firstly let $B_1,\ldots,B_m \subset \mE \times \mE'$ be bBs and
$g:\mR^m \rar \mR$ be a continuous bounded function.
Since $B_i$'s are bounded and $\Phi_i$'s are simple,
given $\Phi_i, i=1,2$, there exists a.s. $N(\Phi_i) \in \mN$ such that
for $n \geq N(\Phi_i)$, the indices $j(n,X_{k}) \neq j(n,X_{l}$)
for $X_{k} \neq X_{l}$, $X_k,X_l \in \Phi_i \cap (B_1 \cup \ldots \cup
B_m)$. Hence for $n\geq N(\Phi_i)$,
$\sE(g(\tilde{\Phi}_i^n(B_1),\ldots,\tilde{\Phi}_i^n(B_m))|\Phi_i) =
\sE(g(\tilde{\Phi}_i(B_1),\ldots,\tilde{\Phi}_i(B_m))|\Phi_i)$ and in
consequence
$\sE(g(\tilde{\Phi}_i^n(B_1),\ldots,\tilde{\Phi}_i^n(B_m))|\Phi_i) \to
\sE(g(\tilde{\Phi}_i(B_1),\ldots,\tilde{\Phi}_i(B_m))|\Phi_i)$ a.s. as
$n \to \infty$.
Since $g$ is bounded, by dominated convergence theorem we have that
$\sE(g(\tilde{\Phi}_i^n(B_1),\allowbreak\ldots,\tilde{\Phi}_i^n(B_m))) \to
\sE(g(\tilde{\Phi}_i(B_1),\ldots,\tilde{\Phi}_i(B_m)))$. Thus
$(\tilde{\Phi}_i^n(B_1),\ldots,\tilde{\Phi}_i^n(B_m)) \cvgdist
(\tilde{\Phi}_i^n(B_1),\ldots,\tilde{\Phi}_i^n(B_m))$. Secondly it is
easy to check that for $B_1 = B' \times B^{''}$,
we have $\sE(\tilde{\Phi}_i^n(B_1)) = \sE(\Phi_i(B'))F(B^{''}) = \sE(\tilde{\Phi}_i(B_1))$ and hence by an
appropriate approximation $\sE(\tilde{\Phi}_i^n(B_1)) =
\sE(\tilde{\Phi}_i(B_1))$ for any bBs $B_1$.

Finally for any bBs $B \subset \mE \times \mE'$
and any realization $\overline Z=\overline z=\{z_{n,j}\}_{n\ge1,j\ge1}$, define
$V_i^{\overline{z}}(B) :=
\int_{\mE}\1[(x,z_{n,j(n,x)}) \in B] \Phi_i(dx)$. Since
$z_{n,j(n,\cdot)}$ is a piecewise constant function,
$1[(x,z_{n,j(n,x)}) \in B]$ is a measurable function in $x$ and so
$V_i^{\overline{z}}$'s are
integral shot-noise fields (as per Definition~\ref{defn:ISN}) indexed by
bBs of $\mE \times \mE'$. Thus from Theorem~\ref{thm:isn_rm}, we have
that for any $dcx$ function $f$,
\begin{eqnarray*}
\lefteqn{\sE(f(\tilde{\Phi}_1^n(B_1),\ldots,\tilde{\Phi}_1^n(B_m))|\overline{Z} = \overline{z}) = \sE(f(V_1^{\overline{z}}(B_1),\ldots,V_1^{\overline{z}}(B_m)))} \\
& \leq &  \sE(f(V_2^{\overline{z}}(B_1),\ldots,V_2^{\overline{z}}(B_m))) = \sE(f(\tilde{\Phi}_2^n(B_1),\ldots,\tilde{\Phi}_2^n(B_m))|\overline{Z} = \overline{z})
\end{eqnarray*}
Now, taking further expectations we get
$(\tilde{\Phi}_1^n(B_1),\ldots,\tilde{\Phi}_1^n(B_m)) \leq_{dcx}
(\tilde{\Phi}_2^n(B_1),\allowbreak\ldots,\tilde{\Phi}_2^n(B_m))$.
Since the approximation satisfies the assumption of Lemma \ref{lem:weak_conv}, the proof follows.
\item  We need to prove
$\sE(f(\overline{\Phi}_1(A_1),\ldots,\overline{\Phi}_1(A_n))) \leq
\sE(f(\overline{\Phi}_2(A_1),\ldots,\overline{\Phi}_1(A_n)))\,$
for $dcx\resp{idc;\,idcv}$ function $f$ and  mutually disjoint $A_k$,
$k=1,\ldots,n$; cf. Proposition~\ref{prop:char_dcx_rm}.
Note that given $\Phi(A_k)=n_k$, we have $\overline{\Phi}(A_k) =
\sum_{i=1}^{n_k}Z^k_i$, where $Z_i^k$ are i.i.d. copies of the
Bernoulli thinning variable. Thus the result follows from the first
statement of Lemma~\ref{lem:LMeester93}.
\item Using the following  fact from
\cite{Muller02}: $X
\leq_{dcx\resp{idcx;\,idcv}} Y$ implies  $X + Z
\leq_{dcx\resp{idcx;\,idcv}} Y + Z$
provided $Z$ is independent of $X$ and  $Y$
one can easily show that  $\Lam_1+\Lam_1' \leq_{dcx\resp{idcx;\,idcv}}
\Lam_2 + \Lam_1'$ assuming  $\Lam_1'$ independent of $\Lam_2$.
The same argument shows that $\Lam_2+\Lam_1' \leq_{dcx\resp{idcx;\,idcv}}
\Lam_2 + \Lam_2'$. The result follows by the transitivity of the order.
\item This result follows easily from Lemma~\ref{lem:weak_conv}
using an increasing approximation of $\mE$ by bBs.
\end{enumerate}
\qed
\end{proof}

\subsection{Impact on Higher Order Properties}
\label{sec:moments}
We will state now some  results involving ordering of
moments of random measures and
draw some conclusions concerning the so called second order
properties. These latter ones make it possible to characterize the clustering
in p.p..

By the $n\,$th power of random measure $\Lambda$, we understand a
random measure $\Lambda^k$ on the product space $\mE^k$
given by $\Lambda^k(A_1\times\ldots\times A_k)=
\prod_{j=1}^k\Lambda(A_j)$.
Its expectation, $\alpha^k(\cdot)=\sE(\Lambda^k(\cdot))$ is called the
$k\,$th moment measure. The first moment measure
$\alpha(\cdot)=\alpha^1(\cdot)$
is called the mean measure.

\begin{prop}\label{prop:product_space}
Consider random measures
 $\Lambda_1\leq_{idcx}\Lambda_2$. Then
$\Lambda_1^k\leq_{idcx}\Lambda_2^k$
and  $\alpha_1^k(\cdot)\le\alpha_2^k(\cdot)$.
 Moreover, if
$\Lambda_1\leq_{dcx}\Lambda_2$ then $\alpha_1(\cdot)=\alpha_2(\cdot)$.
\end{prop}
\begin{proof}
By the standard arguments, one can approximate
any bBs set $C_i$, $i=1,\ldots,n$ in $\mE^k$ by
increasing unions of rectangles. By Lemma~\ref{lem:weak_conv}
and using a similar  argument about composition of a $idcx$ and
$idl$ function as in the proof of Proposition~\ref{prop:char_dcx_rm},
to prove the first statement, it is enough to show the respective
inequality for $idcx$ function
$f:\mR^n\rar\mR$ taken of the values of the moment measures on
$n$ rectangles in $\mE^k$. In this context, consider
$g:\mR^{m}\rar\mR$ given by
$$g(y_1,\ldots,y_m)
=f\Big(\prod_{j\in J_1}y_j,\ldots,\prod_{j\in J_n}y_j\Bigr)\,,$$
where  $J_1,\ldots,J_n$ are $k$-element subsets of
the set $\{1,\ldots,m\}$.
Note for non-negative arguments that  if $f$ is $idcx$  then $g$
is $idcx$.

The second statement follows easily from the first one by the fact
that $f(x)=x$ is $idcx$.
For the first moment (mean measure) note that both $f(x)=x$ and
$f(x)=-x$ are $dcx$.
\qed
\end{proof}

We shall explore now the relation between $dcx$ ordering and
clustering  of points in a p.p.
One of the most popular functions for the analysis of this effect
is  the {\em Ripley's $K$~function} $K(r)$ (reduced second moment function);
see~\cite{Stoyan95}. Assume that
$\Phi$ is a stationary p.p. on $\mR^d$ with finite intensity
$\lambda=\alpha(B)$, where $B$ is a bBs of Lebesgue measure~1.
Then
$$K(r)=\frac1{\lambda|G|}\sE\Bigl(\sum_{X_i\in\Phi\cap G}
(\Phi(B_{X_i}(r))-1)\Bigr)\,,$$
where $B_x(r)$ is the ball centered at $x$ of radius $r$ and
$|G|$ denotes the Lebesgue measure of a bBs $G$;
due to stationarity, the definition does not depend on the choice of $G$.

\begin{prop}\label{prop:Ripley}
Consider two stationary p.p. $\Phi_i$, $i=1,2$,  with same finite intensity
and denote by $K_i(r)$ their Ripley's $K$~functions.
If $\Phi_1\leq_{dcx}\Phi_2$ then $K_1(\cdot)\le K_2(\cdot)$.
\end{prop}
\begin{proof}
Denote $I_i=\sE\Bigl(\sum_{X_j\in\Phi_i\cap G}
(\Phi_i(B_{X_j}(r))-1)\Bigr)$, $i=1,2$.
By the equality of mean measures (Proposition~\ref{prop:product_space}),
it is enough to prove that $I_1\le I_2$.
Note that $I_i$ can be written as the value of some shot noise
evaluated with respect to $\Phi_i^2$, the second product of the p.p..
$$I_i=\sum_{X_j,X_k\in\Phi_i}\1[X_j \in G]\1[0<|X_k-X_j|\le r]\,,$$
where $\1[\cdot]$ denotes the indicator function. Thus, the result follows
from Proposition~\ref{prop:product_space} and Theorem~\ref{thm:isn_rm}.
\end{proof}

Another useful characteristic is the {\em pair correlation
  function} defined on $\mR^2$
as  $g(x,y)=\frac{\rho_2(x,y)}{\rho_1(x)\rho_1(y)},$
where $\rho_k$ is the {\em $k\,$th  product intensity},
equal (outside the diagonals)
to the density of the $k\,$th moment measure $\alpha^k$  with
respect to the Lebesgue measure.

We avoid discussion on questions
such as existence etc.
The following result follows from Proposition~\ref{prop:product_space}.
\begin{cor}\label{cor:pcf}
Consider p.p. such that $\Phi_1\leq_{dcx}\Phi_2$.
Then their respective pair correlation functions satisfy
$g_1(x,y)\le g_2(x,y)$ almost everywhere with respect to the Lebesgue
measure.
\end{cor}

\subsection{Impact on Palm Measures}
\label{sec:order_palm}

For the following definitions and results regarding Palm distributions of
random measures see~\cite[Section~10]{Kallenberg83}.

\begin{defn}
\label{defn:palm_rm}
For a fixed measurable $f$ such that $0 < \sE(\int_{\mE}f(x)\Lam(dx))
< \infty$, the $f$-mixed Palm version of $\Lambda$, denoted by $\Lam_f
\in \mM$, is defined as having the distribution
\[\sP(\Lam_f \in M) = \frac{\sE(\int_{\mE}f(x)\Lam(dx)\1[\Lam \in
    M])}{\sE(\int_{\mE}f(x)\Lam(dx))},\quad M \in \mathcal{M}.\]
In case  $\Lambda$  (say on the Euclidean space  $\mE=\mR^d$)
has a density  $\{\lam(x)\}_{x \in \mR^d}$,
we define for each $x\in\mR^d$
the Palm version $\Lam_x$ of $\Lambda$ by the formula
\[\sP(\Lam_x \in M) 
= \frac{\sE(\lam(x)\1[\Lam \in M])}{\sE(\lam(x))},\quad M \in \mathcal{M}.  \]
\end{defn}
Palm versions  $\Lam_x$ can be defined for a general random
measure via some Radon-Nikodym derivatives.
However, we shall state our result for  $\Lam_x$ as defined above
as well as for  mixed Palm versions  $\Lam_f$
in order to avoid the arbitrariness related to the non-uniqueness
of Radon-Nikodym derivatives.
\begin{prop}
\label{prop:palm_pp_idcx}
Suppose $\Lam_i, i =1,2$ are random measures.
\begin{enumerate}
\item If $\Lam_1\leq_{dcx}\Lam_2$ then
  $(\Lam_1)_f\leq_{idcx}(\Lam_2)_f$
for any non-negative measurable function~$f$ such that $0 <
\int_{\mE}f(x)\alpha(dx)< \infty$, where $\alpha$ is the (common)
mean measure of $\Lambda_i$, $i=1,2$.
\item Suppose that $\Lam_i$ has locally finite mean measure
and  almost surely (a.s.) locally Riemann integrable
density $\lam_i$, $i=1,2$.
If $\{\lam_1(x)\} \leq_{dcx} \{\lam_2(x)\}$, then
$\Lam_1\leq_{dcx} \Lam_2$ and for every  $x \in \mR^d$,
$(\Lam_1)_x \leq_{idcx} (\Lam_{2})_x$.
\end{enumerate}
\end{prop}
\begin{proof}
\begin{enumerate}
\item Denote $I_i=\int_{\mE}f(x)\Lam_i(dx)$, $i=1,2$.
By Proposition~\ref{prop:product_space},
$\Lam_1\leq_{dcx}\Lam_2$
implies that the mean measures are equal and thus
$\sE(I_1)=\sE(I_2)$.
It remains to prove
$$\sE(g(\Lambda_1(B_1),\ldots,\Lambda_1(B_n))I_1)
\le\sE(g(\Lambda_2(B_1),\ldots,\Lambda_2(B_n))I_2)$$
for $idcx$ function $g$.
This follows from Theorem~\ref{thm:isn_rm}
and the fact that
$h(x_0,x)= x_0g(x):\mR^{n+1} \to \mR$
is $idcx$, for non-negative argument $x_0$.
\item The first part follows immediately from the second statement of
Lemma~\ref{lem:Miyoshi04a}.
For the second part, use the same argument about $h(x_0,x)= x_0g(x)$ as
above.
\remove{here $a^+$ for $a \in \mR$ denotes the non-negative part of $a$. Check that if $f$ is twice differentiable,
then so is $g$ for $x_0 \geq 0$ and all the second partial derivatives of $g$ are positive. Thus by Theorem  from \cite{Muller02},
$dcx$ condition is satisfied for $x_0 \geq 0$. For $x_0 < 0$, the $dcx$ condition can be checked directly. Since twice differentiable
$dcx$ functions generate the $dcx$ order, this completes the proof.}
\end{enumerate}

\qed
\end{proof}

\begin{rem}\label{rem:Palm_counterexample}
Compared to earlier results where $dcx$ ordering led to $dcx$ ordering, one might tend to believe
that the loss here (as $dcx$  implies $idcx$ only)
is more technical. However the following
illustrates that it is natural to expect so: consider a Poisson
p.p. $\Phi$ and its (deterministic) intensity measure $\alpha(\cdot)$
(i.e., its mean measure $\alpha(\cdot)=\sE(\Phi(\cdot))$.
Using the complete independence property of the Poisson p.p.
and the fact that each $dcx$ function is component-wise convex,
one can show that for  disjoint bBs $A_1,\ldots,A_n$ and any $dcx$
function $f$,
$f(\alpha(A_1),\ldots,\alpha(A_n))\le \sE(f(\Phi(A_1),\ldots,\Phi(A_n))$.
Thus $\alpha\leq_{dcx}\Phi$. It is easy to see that
$\alpha_f(\cdot)=\alpha(\cdot)$ (mixed Palm version of a
deterministic measure is equal to the original measure).
Take $f(x)=1[x\in A]$ for some bBs $A$.
Then $\sE(\Phi_f(A))=\sE((\Phi(A))^2)/\alpha(A)=\alpha(A)+1$
since $\Phi(A)$ is a Poisson r.v..
Thus $\alpha_f(A)<\sE(\Phi_f(A))$
disproving   $\alpha_f(A)\leq_{dcx} \Phi_f(A)$.
Another counterexample involving Poisson-Poisson  cluster p.p. will
be given in Remark~\ref{rem:CoxPalm}.
\end{rem}

\subsection{Cox Point Processes}
\label{ss.Cox}
We will consider now Cox p.p. (see e.g.~\cite[III
5.2]{Stoyan95}),
known also as doubly stochastic Poisson p.p., which constitute
a rich class  often used to model patterns which
exhibit more clustering than in Poisson p.p..

Recall that a {\em Cox$\,(\Lambda)$ p.p.}
 $\Phi_{\Lambda}$ on $\mE$ generated
by the random {\em intensity measure} $\Lambda(\cdot)$ on $\mE$
is defined as having the property that
$\Phi_{\Lambda}$ conditioned on $\Lambda(\cdot)$
is a Poisson p.p. with intensity
$\Lambda(\cdot)$. Note that  Cox p.p. may be
seen as a result of an operation transforming some  random (intensity) measure
into a point (Cox) p.p..

One can easily show that this operation preserves our orders.
\begin{prop}
\label{prop:int_fld_meas}
Consider two ordered random measures $\Lambda_1 \leq_{dcx\resp{idcx;\,idcv}}
\Lambda_2$. Then $\Phi_{\Lambda_1} \leq_{dcx\resp{idcx;\,idcv}} \Phi_{\Lambda_2}$.
\end{prop}
\begin{proof}
Taking a $dcx\resp{idcx;\,idcv}$ function $\phi$,
assuming (by Proposition~\ref{prop:char_dcx_rm})
mutually disjoint bBs $A_k$, $k=1,\ldots,n$, using the
definition of  Cox p.p. and the
second statement of the Lemma~\ref{lem:LMeester93}
one shows for $i=1,2$ that  that
the conditional expectation
$$\sE(\phi(\Phi_{\Lambda_i}(A_1),\ldots,\Phi_{\Lambda_i}(A_n))|\Lambda_i)$$
given the intensity measure $\Lambda_i$ is a  $dcx\resp{idcx;\,idcv}$
function of
$(\Lambda_i(A_1),\ldots,\Lambda_i(A_n))$.
The result follows thus from the assumption of the
measures $\Lambda_i$ being $dcx$ ordered.
\qed
\end{proof}

We will show now using Theorem~\ref{thm:isn_rm}
that $dcx, idcx,idcv$ ordering of Cox intensity measures
is preserved by independent (not necessarily identically
distributed) marking and thinning, as well as
independent displacement of points of the p.p..


By {\em independent  marking} of p.p. $\Phi$ on $\mE$
with {\em marks}  on some LCSC space $\mE'$,
we understand a p.p. $\tilde{\Phi} =
\sum_{i}\varepsilon_{(x_i,Z_i)}$
such  that
given $\Phi=\sum_{i}\varepsilon_{x_i}$,
$Z_i$ are independent random elements in $\mE'$,
with distribution $\sP\{Z_i\in\cdot|\Phi= \sum_{i}\varepsilon_{x_i}\}
=F_{x_i}(\cdot)$ given by some probability (mark) kernel $F_x(\cdot)$ from
$\mE$ to $\mE'$. The fact that  $F_x(\cdot)$ may depend on $x$
(in contrast to i.i.d. marking) is
sometimes emphasized by calling  $\tilde{\Phi}$ a ``position
dependent'' marking.
Independent 
thinning can be
seen as the  projection on $\mE$ of the subset $\tilde\Phi(\cdot,\{1\})$ of
the independently marked  p.p. $\tilde\Phi$
where the marks $Z_i\in\{0,1\}=\mE'$, are independent  Bernoulli
thinning variables $Z_i=Z_i(x)$, whose distributions may be dependent
on $x_i$.
Similarly, the projection of an independently marked p.p. $\tilde\Phi=\sum_i
\varepsilon_{(x_i,Z_i)}$ on the space of
marks $\mE'$; i.e., $\tilde\Phi(\mE\times\cdot)=\sum_i
\varepsilon_{Z_i}$ can be seen as independent displacement of points of $\Phi$
to the space $\mE'$. Special examples are i.i.d. shifts of points in
the Euclidean space, when $Z_i=x_i+Y_i$, where $Y_i$ are i.i.d.

\begin{prop}\label{prop:mark_cox}
Suppose $\Phi_i, i =1,2$, are two Cox$\,(\Lam_i)$   p.p..
Assume that their intensity measures are ordered $\Lam_1
\leq_{dcx\resp{idcx;\,idcv}} \Lam_2$.
Let $\tilde{\Phi}_i, i =1,2$ be the corresponding independently
marked p.p. with the same mark kernel  $F_x(\cdot)$.
Then $\tilde{\Phi}_1 \leq_{dcx\resp{idcx;\,idcv}} \tilde{\Phi}_2$.
\end{prop}

{}From the above Proposition,  the following corollary follows immediately
by the last  statement of Proposition~\ref{prop:mpp_idcx}.
\begin{cor}
Independent thinning and  displacement of points
preserves $dcx\resp{\allowbreak idcx;\,idcv}$ order of the intensities of Cox p.p..
\end{cor}

\begin{proof} (\emph{Prop.~\ref{prop:mark_cox}})
Let $\Phi_i$ be Cox$\,(\Lam_i)$ $i=1,2$ respectively.
 Assume
 $\Lam_1 \leq_{dcx(idcx,idcv)} \Lam_2$.
It is known that independent marking of Cox$\,(\Lam_i)$
p.p. is a Cox$\,(\tilde\Lam_i)$ p.p.
with intensity measure $\tilde\Lam_i$ on $\mE\times\mE'$
given by $\tilde\Lam_i(\cdot)=\int_{\mE}\int_{\mE'}\1[(x,y)\in\cdot]\,F_x(d
y)\Lam_i(d x)$; cf.~\cite[Secs 4.2 and 5.2]{Stoyan95}.
Let $S$ be the family of bBs in $\mE\times\mE'$;
for  $x\in\mE$ and bBs $C\subset\mE\times\mE'$
consider $h(x,C)=\int_{\mE'}\1[(x,y)\in
C]\,F_x(d y)$.
Then the integral shot noise $V_{\Lam_i}(C)=\int_{\mE}h(x,C)\,\Lam_i(d
x)$ satisfies $V_{\Lam_i}(C)=\tilde\Lam_i(C)$ for all bBs $C$.
Thus, by Theorem~\ref{thm:isn_rm}
$\tilde\Lam_1\leq_{dxc\resp{idcx;\,idxv}}\tilde\Lam_2$
and the result follows from Proposition~\ref{prop:int_fld_meas}.
\qed
\end{proof}

\remove{We shall now state a corollary that is very useful for comparing
many useful quantities.
\begin{cor}
\label{cor:intl_fnl_isn} Let the assumptions be as in Theorem
\ref{thm:isn_rm} with $S = \mR^d$. If $\{f_y\}_{y \in \mR^d}$ be
family of $icx$ functions from $\mR$ to $\mR$, then $\{f_y(V^1(y))\
\leq_{idcx} \{f_s(V^2(s))\}.$ And if the two random fields are a.s.
Riemann integrable
\[
\left( \int_{I_1}f_y(V^1(y))dy,\ldots,\int_{I_n}f_y(V^1(y))dy
\right) \leq_{idcx}
\left(\int_{I_1}f_y(V^2(y))dy,\ldots,\int_{I_n}f_y(V^2(y))dy \right)
\]
for $I_1,\ldots,I_n$ disjoint bBs in  $\mR^d$.
\end{cor}

The first result follows due to composition rules of idcx functions
and the second result follows directly from Lemma \ref{lem:Miyoshi04a}.

Now we state the ordering of shot-noise fields results for the Palm version which follows
from Proposition \ref{prop:palm_pp_idcx}.

\begin{cor}
\label{cor:palm_isn}
Let $\Phi_i,i=1,2$ be two Cox($\Lam_i$).
\begin{enumerate}
\item  If $\Lam_1 \leq_{dcx} \Lam_2 $, then $\{V_{\Lambda_{1f}}(x)\} \leq_{idcx} \{V_{\Lambda_{2f}}(x)\}$ and
$\{V_{\Phi_{1f}}(x)\} \leq_{idcx} \{V_{\Phi_{2f}}(x)\}$ .
\item Suppose that $\{\lam_1(x)\} \leq_{dcx} \{\lam_2(x)\}$ are the two stationary density fields of $\Lam_i$ respectively.
Then for every $s \in \mR^d$, $\{V_{\Lambda_{1s}}(x)\} \leq_{idcx} \{V_{\Lambda_{2s}}(x)\}$ and $\{V_{\Phi_{1s}}(x)\} \leq_{dcx} \{V_{\Phi_{2s}}(x)\}.$
\end{enumerate}
\end{cor}
}

If $\Lambda(\cdot) \in \mM(\mR^d)$ a.s has a density $\{\lambda(x)\}_{x
\in \mR^d}$ with respect to Lebesgue measure then the density is
referred to as the intensity field of the Cox p.p., which will be
called in this case Cox$\,(\lam)$  p.p. and denoted by $\Phi_\lam$.

It is known that Cox p.p. is  over-dispersed with respect to
the Poisson p.p., i.e., $\text{Var}(\Phi_1(B))\le\text{Var}(\Phi_2(B))$
where $\Phi_1,\Phi_2$ are, respectively, Poisson and Cox p.p. with the
same mean measure. Hence, it is clear that a Cox p.p. can only be greater
in $dcx$ order than a Poisson p.p. with the same mean measure.
Indeed, in Section~\ref{sec:examples}  we will show several examples when
this stronger result holds, namely  Cox
p.p. that are $dcx$ ordered (larger) with respect to the
corresponding Poisson p.p., as well as Cox p.p. $dcx$ ordered with
respect to each other.

\remove{
\begin{prop}
\label{prop:branching_cox} Consider two Euclidean p.p.
 $\Phi_1 \leq_{dcx(idcx,idcv)} \Phi_2$
Euclidean p.p.
and $\eta$ is another Euclidean p.p. independent of the above. Let
$\widehat{\Phi}_i, i =1,2$ be
p.p. defined as follows $\widehat{\Phi}_i = \cup_{X_j \in
\Phi_i} (X_j + \eta_j)$ for $i=1,2$ where $\eta_j$ are i.i.d copies
of $\eta$. We further assume one of the following statements is
true.

\begin{itemize}
\item $\Phi_i(\mR^d) < \infty$ a.s. for $i=1,2$ and $\eta$ is a
stationary p.p.

\item $\Phi_i, i=1,2$ are Cox p.p..

\end{itemize}
Then $\Phi_1 \leq_{dcx(idcx,idcv)} \Phi_2$ implies that
$\widehat{\Phi}_1 \leq_{dcx(idcx,idcv)} \widehat{\Phi}_2$

\end{prop}

\begin{proof} Firstly note that the assumptions imply that the new
p.p. is well-defined. Observe that under the first
assumption for $i=1,2$, $\widehat{\Phi}_i = \cup_{X_j \in \Phi_i}
(X_j + \eta_j) =^d  \cup_{X_j \in \Phi_i} \eta_j,$ where the last
inequality is in distribution. Hence for a bBs $B$
and $i=1,2$, we get $\widehat{\Phi}_i(B) = \sum_{j=1}^{\Phi_i(\mR^d)} \eta_j(B).$
Now the proof follows similar to the above two propositions.

Lets now prove under the second assumption. Pick a sequence of
compact sets $A_k, k \in \mN$ such that $A_k \nearrow \mR^d.$ For
$i=1,2$ and $k \in \mN$, define $\widehat{\Phi}^k_i(B) = \sum_{X_j
\in \Phi_i \cap A_k} (\eta_j-X_j)(B)$. For $i=1,2$ note that
conditioned on $\Phi_i(A_k)$ , $(\eta_j-X_j)_{1 \leq j \leq
\Phi_i(A_k)}$ are i.i.d point p.p.. .............. \qed
\end{proof}
}

\subsection{Alternative Definition of $dcx$ Order}
\label{sec:alt_defn}

We viewed a random measure as a random field and have defined
ordering from this viewpoint.
Alternatively, one can consider a
random measure as an element of the space of Radon measures $\mM$
and define ordering between two
$\mM$-valued random elements. This can be done once we define what
is a $dcx$ function on $\mM$. The $dcx$ order can be defined on
more general spaces; \cite{LMeester99} extends the notion of $dcx$
ordering to lattice ordered Abelian semigroups
with some compatibility conditions between the lattice structure and
the Abelian structure ($LOAS^+$).
The space $\mM$ can be equipped with the
following lattice and algebraic structure. Consider  the following partial order:
for $\mu,\nu \in \mM$, we say $\mu \leq \nu$ if $\mu(B) \leq
  \nu(B)$ for all bBs $B$ in $\mE$ and addition
$(\mu + \nu)(B) = \mu(B) + \nu(B)$.
Under this definition, the space
$\mM$ forms a $LOAS^+$ as required by \cite{LMeester99}. Then
one can define a directionally convex function on $\mM$ as in
Definition~\ref{defn:dcx_fn}. Call it a $dcx^1$ function.
This gives rise to  $dcx^1$ order of random measures analogously to
the first part of the Definition~\ref{defn:dcx_order}.

\remove{
\begin{ex}
\label{ex:defn2} Let $\{Y_i\}_{i \geq 1}$ be i.i.d sequence of
random variables such that $\sP(Y_i=0) = \sP(Y_i=2) = 1/2.$
Define two p.p. $\Phi_1$ and $\Phi_2$ on $\mR$
as follows $\Phi_1=\sum_{i=1}^\infty\varepsilon_{i}$,
$\Phi_2=\sum_{i=1}^\infty Y_i\varepsilon_{i}$.
 As a simple consequence of Jensen's inequality, it
follows that $\Phi_1 \leq_{dcx}\Phi_2$. Define the
following functional of the p.p. (random measure)
$g(\Phi) =\liminf_{n\to\infty}\Phi(\{n\})$.
\selfnote{Is this functional $dcx^1$?}
Now see that $\sE(g(\Phi_1))= 1 > 0 = \sE(g(\Phi_2))$.
 \end{ex}
}

\remove{The converse isn't clear. This is an incarnation
of a more general question:
For any discrete time stochastic process, does the $dcx$ ordering
as a stochastic imply $dcx$ ordering as an infinite sequence?
\cite{Kamae77} proves so for increasing order and this is used in a
method very similar to one sketched below to prove the equivalence
of two definitions for strong ordering by \cite{Rolski91}.}

Now we have two reasonable definitions of ordering of random
measures. It is easy to see that $dcx^1$ ordering implies $dcx$
ordering. In light of Example 5.1.7 of \cite{Muller02}, existence of a counterexample
to the converse looks plausible, though we failed in our attempts to construct one.
However, the result of \cite{Bassan91} proves that {\em convex} ordering of
real valued stochastic process $\{X_n\}_{n\in\mN}$
implies {\em continuous, convex}
ordering of the corresponding elements of the
infinite-dimensional Euclidean spaces $\mR^\mN$.
This suggests  that
$dcx$  of random measures may imply  a $dcx^{1*}$ order
induced by some subclass  of $dcx^1$ functionals of random measures,
 which are regular in some sense. Leaving this general question
as an open problem, we remark only that
the integral shot-noise fields  studied in the next
section can be seen as some particular class of
functionals of random measures, which are $dcx^1$ (in fact linear on
$\mM$) and regular enough for their
means to satisfy the required inequality
provided the random measures are $dcx$ ordered.
It is natural thus to have them in the suggested $dcx^{1*}$ class.

\remove{For the
specific problem of random measures, the following method adapted
from proof of \cite{Rolski91} does half the job. The sketch of
possible proof is as follows : Let $\texttt{I} = (I_1,\ldots)$ be
the countable DC-semiring of $\mE$ (See \cite{Kallenberg83}.)
Assume that $\Lambda_1 \leq_{dcx} \Lambda_2$ implies
$\Upsilon(\Lambda_1) \leq_{dcx^2} \Upsilon(\Lambda_2)$ where
$\Upsilon(\mu) = (\mu(I_1),\ldots)$ is the infinite sequence and
$dcx^2$ stands for $dcx$ ordering as an element of $R^{\infty}_+$.
It is easy to show that for $f:\mM \rar R$ $dcx^1$, $f.\Upsilon^{-1}
:\Upsilon(\mM) \rar R$ is $dcx^2$. Hence to prove $dcx$ ordering
implies $dcx^1$ ordering, it suffices to prove that $dcx$ ordering
implies $dcx^2$ ordering.}

Recall also that for {\em strong} order of p.p.
there is the full equivalence between these two definitions,
and both imply the possibility of a coupling of the ordered p.p.
such that the smaller one  is a.s. a subset of the greater one;
cf~\cite{Rolski91}.

\section{Ordering of Shot-Noise Fields}
\label{sec:OSN}
In this section we will prove Theorem~\ref{thm:isn_rm}
concerning $dcx$ ordering of integral shot-noise fields, which is the main
result of this paper. We will also consider the so called extremal
shot-noise fields.
\subsection{Integral Shot-Noise Fields}
\label{sec:ISN}
Usually shot-noise fields are defined for p.p.
as the following sum (thus sometimes called additive shot-noise fields)
$ V_{\Phi}(y) = \sum_{X_n \in \Phi}h(X_n,y)$
where $\Phi=\sum_n\varepsilon_{X_n}$
and $h$ is a non-negative response function. In definition~\ref{defn:ISN}
we have made a significant but natural generalization of this definition.
It is pretty clear as to why we call this generalization
integral shot-noise field.
The extension to unbounded response functions is not just a mathematical generalization alone. It shall
provide us a simple proof of ordering for extremal-shot-noise fields for p.p..

Now, we shall  prove Theorem~\ref{thm:isn_rm}.
The proof is inspired by \cite{Miyoshi04}.
\begin{proof}(\emph{Theorem \ref{thm:isn_rm}})
We shall prove the second statement first. The necessary
modifications for the proof of the first statement  shall be indicated
later on.
\begin{itemize}
\item[{\rm 2.}] We need to show that $(V^1(y_1),\ldots,V^1(y_m)) \leq_{dcx} (V^2(y_1),\ldots,V^2(y_m))$
for $y_i \in S, 1 \leq i \leq m$ and $V^j(\cdot) = V_{\Lambda_j}(\cdot)$,
$j=1,2$. The proof relies on the construction of  two sequences of random
vectors $(V^j_k(y_1),\ldots,V^j_k(y_m))$, $k=1,2.\ldots$, $j=1,2$ satisfying the
assumptions of Lemma \ref{lem:weak_conv}.

Choose an increasing sequence of compact sets $K_k$, $k \geq 1$ in
$\mE$,  such that
$K_k \nearrow \mE$. Since $h$ is measurable in its first argument,
we know that there exists a sequence of
simple functions $h_k(\cdot,y_i), k \in \mN$ such that as $k \rar
\infty$, $h_k(\cdot,y_i) \uparrow h(\cdot,y_i)$ for $1\leq i \leq m$. They
can be written down explicitly as follows:
\begin{eqnarray*}
h_k(\cdot,y_i)
&=&
\gamma_k\1[\{x\in K_k:
  h(x,y_i)=\infty\}]\\[1ex]
&&+\sum_{n=1}^{\gamma_k} \frac{n-1}{2^k}
\1[\{x\in K_k:\frac{n-1}{2^k} \leq h(x,y_i) <
    \frac{n}{2^k}\}](\cdot)
\end{eqnarray*}
for $1 \leq i \leq m$, where $\gamma_k = k2^k$.
Put $I^i_{kn} = \{ x \in K_k :
\frac{n-1}{2^k} \leq h(x,y_i) < \frac{n}{2^k} \}$
and $I^i_{k\infty} =  \{ x \in K_k: h(x,y_i) =\infty \}$
for $1 \leq i \leq
m$ and $1 \leq n \leq \gamma_k$. Note that all $I^i_{kn}$
$n=1,\ldots,\infty$ are bBs
and the  sequence of random vectors we are looking for is
\[ V^j_k(y_i) = \int_{\mE}h_k(x,y_i) \Lambda_j(dx) =
\gamma_k\Lambda_j(I^i_{k\infty})+ \sum_{n=1}^{\gamma_k} \frac{n-1}{2^k}
\Lambda_j(I^i_{kn}), \]
for $j=1,2$.
By the definition of integral, it is clear that
for $j=1,2$ as $k \rar \infty$, $(V^j_k(y_1),\ldots,V^j_k(y_m))
\uparrow (V^j(y_1),\ldots,V^j(y_m))$ a.s. and hence in distribution.
By monotone convergence theorem, the expectations, which are finite by
the assumption, also converge.
What remains to prove is that for each $k \in \mN$, the vectors are
$dcx$  ordered.

Fix $k \in \mN$.  Now observe that for $j=1,2$, $i=1,\ldots,m$,
$V^j_k(y_i)$ are increasing linear functions of the vectors $(\Lambda_j(I^i_{kn}):
n=1,\ldots,\gamma_k,\infty)$, $j=1,2$.
The latter are $dcx$ ordered by the assumptions. And
since composition of $dcx$ with increasing linear
functions is $dcx$,  it follows that $(V^1_k(y_1),\ldots,V^1_k(y_m)) \leq_{dcx}
(V^2_k(y_1),\ldots,V^2_k(y_m)).$\\

\item[{\rm 1.}] For vectors $(V^j_k(y_1),\ldots,V^j_k(y_m))$, $k=1,2.\ldots$, $j=1,2$ defined as above,
$f(V^j_k(y_1),\ldots,V^j_k(y_m)) \uparrow f(V^j(y_1),\ldots,V^j(y_m))$ a.s. for $f$ $idcx\resp{idcv}$
and hence \\
$\sE(f(V^j_k(y_1),\ldots,V^j_k(y_m))) \uparrow \sE(f(V^j(y_1),\ldots,V^j(y_m)))$, $j=1,2$.
The proof is complete by noting that $\sE{f(V^1_k(y_1),\ldots,V^1_k(y_m))} \leq
\sE{(f(V^2_k(y_1),\ldots,V^2_k(y_m))}$ for all $k \geq1$ and $f$
$idcx\resp{idcv}$.
\end{itemize}
\qed
\end{proof}

\subsection{Extremal  Shot-Noise Fields}
\label{sec:Max_SN}

We recall now the definition of the extremal shot-noise, first introduced in \cite{Heinrich94}.
\begin{defn}\label{defn:ESN}
Let $S$ be any set and $\mE$ a LCSC space. Given a p.p.
$\Phi$ on $\mE$ and a measurable (in the first variable alone)
response function $h(x,y) :\mE \times S \rar \mathbb{R}$,
the extremal shot-noise field is defined as
\begin{equation}
\label{eqn:esn} U_{\Phi}(y) = \sup_{X_i \in \Phi}\{h(X_i,y)\}.
\end{equation}
\end{defn}
In order to state our result for extremal shot-noise fields, we
shall use the {\em lower orthant}~($lo$) order.
\begin{defn}
\label{defn:lo}
Let $X$ and $Y$ be random $\mR^d$ vectors. We say $X \leq_{lo} Y$ if
$\sP(X \leq t) \geq \sP(Y \leq t)$ for every $t \in \mR^d$.
\end{defn}
On the real line, this is the same as strong order (i.e., when $\mathfrak{F}$
consists of increasing functions) but in higher dimensions it is
different. Obviously  $st$ order
implies $lo$ order and examples of  random vectors which are
ordered in $lo$ but not in $st$ are known; see (\cite{Muller02}).
Thus, it is clear that the following proposition is a generalization of
the corresponding one-dimensional result in \cite{Miyoshi04} where the proof
method was similar to the proof of the ordering of integral shot-noise fields.
We shall give a much simpler proof using the already proved result.
\begin{prop}
\label{prop:lo_msn} Let $\Phi_1 \leq_{idcv}
\Phi_2$. Then $\{U_{\Phi_1}(y)\}_{y \in S} \leq_{lo}\{U_{\Phi_2}(y)\}_{y \in S}$.
\end{prop}
\begin{proof}
The probability distribution function of the
extremal shot-noise can be expressed
by the Laplace transform of some corresponding (additive) one
as follows.
Let $\{x_1,\ldots,x_m\} \subset S$ and
  $(a_1,\ldots,a_m) \in R^m$. Then
\begin{eqnarray*}
\sP(U(y_i) \leq a_i , 1 \leq i \leq m) & = & \sE(\prod_i
\1[\sup_n\{h(X_n,y_i) \leq a_i\}]) \\
& = & \sE(\prod_i \prod_n \1[h(X_n,y_i) \leq a_i]) \\
& = & \sE(\prod_i \prod_n e^{\log{\1[h(X_n,y_i) \leq a_i]}}) \\
& = & \sE(\prod_i  e^{- \sum_n-\log{\1[h(X_n,y_i) \leq a_i]}}) \\
& = & \sE( e^{- \sum_i \hat{U}(y_i)})
\end{eqnarray*}
where $\hat{U}(y_i) = \sum_n-\log{\1[h(X_n,y_i) \leq a_i]}$ is an
additive shot-noise with response function taking values in
$[0,\infty].$ The response function is clearly non-negative and
measurable. The function $f(x_1,\ldots,x_m) = e^{-\sum_i
x_i}$ is a ddcx function on $(-\infty,\infty]$.
The result follows by the first statement of
Theorem~\ref{thm:isn_rm}.
\qed
\end{proof}

The extremal shot-noise field can be used to define the Boolean model.
Given a (generic) random closed set (RACS; see~\cite[Ch.~6]{Stoyan95}) $G$,
let $h((x,G),y) = \1[y \in x + G]$.
\begin{defn}
\label{defn:boolean_model}
By a {\em Boolean model} with the p.p. of germs $\Phi$  and
the typical grain $G$ we call  the random set
$C(\Phi,G) = \{y: U_{\tilde{\Phi}}(y) > 0\}$~%
where
$\tilde{\Phi} = \sum_i \varepsilon_{(X_i,G_i)}$ is i.i.d. marking of
$\Phi$ with the mark distribution equal to this of $G$.
\end{defn}
We shall call $G$ a fixed grain if there exists a closed set $B$ such
that $G = B$ a.s.. We shall demonstrate in Section \ref{ssec:coverage}
as to how one can obtain comparison results for the Boolean model
using the results of this section.

\section{Examples of $dcx$ Ordered Measures and Point Processes}

\label{sec:examples}

In this section, we shall provide some examples of $dcx$ ordered
measures and p.p. on the Euclidean space
$\mE = \mR^d$.  The examples are intended to be illustrative and not
encyclopaedic. The purpose of the examples is to show that
there are $dcx$ ordered p.p. as well as demonstrate some methods
to prove that two p.p. are $dcx$ ordered. Many of the examples seem to
indicate that p.p. higher in $dcx$ order cluster more, at
least for Cox p.p..
\remove{Also, under suitable assumptions one can
obtain a class   $dcx$ or $idcx$ ordered p.p.. There are various
other results (for e.g. ordering of Palm
versions of the p.p.) possible combining the results on
$dcx$ order and the results in our previous sections. However, we
refrain from mentioning them as they are fairly obvious.}

\subsection{Ising-Poisson Cluster Point Processes}
\label{sec:ising_cox}

Let $\{\lambda(s)\}_{s \in \mR^d}$ be a stationary random intensity
field. Define a new field, which is random but constant in space
$\{\lambda_m(s)=\lambda(0)\}$ and deterministic constant field
$\{\lambda_h(s) = \sE(\lambda(0))\}$.
Cox($\lambda_m$) is known as mixed
Poisson p.p. and Cox($\lambda_h$) is just the well-known homogeneous Poisson
p.p.. Denote the random intensity measures of the Cox, mixed and homogeneous
Poisson p.p.,  by $\Lambda,\Lambda_m$ and $\Lambda_h$ respectively
(i.e., $\Lambda(dx)=\lambda(x)\,dx$, etc.)
It is proved in~\cite{Miyoshi04} that $\Lambda   \leq_{dcx} \Lambda_m $ and
when $\{\lambda(s)\}$ is a conditionally increasing field, $\Lambda_h
\leq_{dcx}\Lambda$. Recall that a random field $\{X(s)\}$ is a
{\em conditionally increasing field} if
for any $k$ and $s_1,\ldots,s_k\in\mR^d$ the expectation
$\sE(f(X(s_1))|X(s_j) = a_j \ \forall \ 2 \leq j \leq k)$ is increasing in $a_j$
for all increasing $f$.
However, no example of a conditionally increasing field was given
in~\cite{Miyoshi04}. Now we construct one.

Consider the $d$-dimensional lattice $\mZ^d$. Let
$\{X(z)\}_{z\in\mZ^d}$ be i.i.d. random variables taking values
in~$\{+1,-1\}$. Call $\{X(z)\}$ a (random) configuration of spins.
In order to obtain a stationary field  consider a random shift
of  the origin of $\mZ^d$ to $U$
with uniform distribution on $[0,1]^d$ ($U$ independent of
$\{X(z)\}$).  Let the lattice shifted by
$U$ be denoted by $\mZ^d_*$. Pick two numbers $\mu_2 \leq \mu_1.$
For $s\in\mR^d$, define $\lambda(s) = \mu_1 \1[X(\dot{s})=1]+ \mu_2
1[X(\dot{s})=-1]$ where $\dot{s}$ represents the unique
``lower left'' point in $\mZ^d_*$ nearest to $s$.
The intensity field is clearly stationary.
We shall now show that $\{\lam(s)\}$ is
conditionally increasing.
Note that
\begin{equation}
\label{eqn:ising_intensity}
f(\lambda(s))  =
1[x(\dot{s})=1](f(\mu_1)-f(\mu_2)) + f(\mu_2)
\end{equation}
{}From Theorem 1.2.15 of \cite{Muller02}, it is sufficient to show the
conditional increasing property conditioned
on $U$, the random origin of the lattice $\mZ^d_*$.
Hence it is enough for the
Ising model to possess the following property:
\begin{eqnarray*}
\lefteqn{ \sP(X(z_1)=1 | X(z_2) = -1,X(z_j)=a_j, j=3,\ldots,k)}\\
&\leq &\sP(X(z_1)=1 | X(z_2) = 1,X(z_j)=a_j, j=3,\ldots,k),
\end{eqnarray*}
where $a_i \in \{+1,-1\}$ and $z_i \in \mZ^d, i=1,\ldots,k$.
This  follows easily from the fact that the spins are i.i.d.

We call the Cox  p.p.
generated by the above conditionally increasing field $\{\lambda(s)\}$
the {\em Ising-Poisson cluster} p.p. By the arguments
presented in~\cite{Miyoshi04},  it is $dcx$ larger than the
homogeneous Poisson p.p. with the same intensity.
Note that intuitively the Ising-Poisson cluster  p.p.
``clusters'' its points  more than a homogeneous Poisson p.p.
\remove{A slightly specialized version of
this shall illustrate the point better. Put $\mu_2 = 0$. Then
$\bar{\lambda} = \mu_1p$ and if $p < 1$, $\bar{\lambda} < \mu_1$.
Then the Ising-Poisson cluster
p.p. shall have high intensity in boxes where
the lower left vertex has a positive spin and no point otherwise.
And thus one would see the points in clusters.}
In what follows, we will see more examples of cluster (Cox) p.p. which are $dcx$
larger than the corresponding homogeneous Poisson p.p..

\subsection{L\'{e}vy Based Cox Point Processes (LCPs)}
\label{sec:lcp}

This class of p.p. is being introduced in
\cite{Hellmund08}. One can find many examples of LCPs in the above
mentioned paper. In simple terms, a LCP is a p.p. whose
intensity field is an integral shot-noise field of a L\'{e}vy basis. A random
measure $L \in \mM(\mR^d)$ is said to be a non-negative L\'{e}vy
basis if
\begin{itemize}

\item for any sequence $\{A_n\} $ of disjoint, bBs of $\mR^d$,
$L(A_n)$ are independent random variables ({\em complete
  independence\,}) and $L(\bigcup A_n) = \sum L(A_n)$ a.s. provided
$\cup A_n$ is also a bBs of $\mR^d$.

\item for every bBs $A$ of $\mR^d$, $L(A)$ is infinitely divisible.

\end{itemize}
We shall consider only non-negative L\'{e}vy bases, even though there exist signed L\'{e}vy bases too
(see \cite{Hellmund08}). Hence, we shall omit the reference to
non-negativity in future.

A Cox p.p. $\Phi$
is said to be a LCP, if its intensity field is of the form
\[ \lambda(y) = \int_{\mR^d}k(x,y)L(dx), \]
where $L$ is a L\'{e}vy basis and the kernel $k$ is a non-negative function such that
$k(x,y)$ is a.s. integrable with respect to $L$ and $k(.,y)$ is
integrable with respect to
Lebesgue measure.
In~\cite{Hellmund08} the response function $k$ and the L\'evy basis
$L$ is chosen such that $\int_B\lambda(y)\,d y<\infty$ a.s. for all bBs
$B$, for which a sufficient condition is  $\int_B\sE(\lambda(y))\,d
y<\infty$. In our considerations, in order to be able to use
Lemma~\ref{lem:Miyoshi04a}, we will require that $\lambda(y)$ is
a.s. locally Riemann integrable.

\begin{rem}
Note that a sufficient condition
for this is that $\lambda(y)$ is a.s. continuous,
for which, in turn, it is enough  to assume
that  $k$ is continuous in its second argument
and that for all $x\in\mR^d$,  there exist
$B_{x}(\epsilon_x)$, $\epsilon_x>0$ such that
$\int_{\mR^d}\sup_{z\in B_{x}(\epsilon_x)}k(z,y)\,\alpha(dx)
\allowbreak<\infty$
for all $y$, where  $\alpha(B) = \sE(L(B))$,
the mean measures of the L\'{e}vy bases;
(cf~\cite{Baccelli01}).
\end{rem}

\begin{lem}
\label{lem:lbasis}
Let $L_1$ and $L_2$ be L\'{e}vy bases with mean measure $\alpha_i$.
Let $\Phi_i, i=1,2$ be LCPs with
L\'{e}vy bases $L_i, i=1,2$ respectively.
\begin{enumerate}
\item $L_1 \leq_{dcx\resp{idcx;\,idcv}} L_2$ if and only if $L_1(A)
\leq_{cx\resp{icx;\,icv}} L_2(A)$
for all bBs $A$ of $\mR^d$,  where $cx, icx,icv$ stands, respectively
for convex, increasing convex and  increasing concave.

\item If $L_1 \leq_{dcx\resp{idcx;\,idcv}} L_2$, then  $\Phi_1
  \leq_{dcx\resp{idcx;\,idcv}} \Phi_2$ provided  the intensity fields
  $\lambda_i(y)$ of LCP $\Phi_i$
is a.s. locally Riemann integrable
with  these   integrals,  in case of $dcx$,  having finite means.

\item $\alpha_i \leq_{dcx} L_i$.
\end{enumerate}

\end{lem}
\begin{proof}
The first part is due to  Proposition~\ref{prop:char_dcx_rm} and
the complete independence property of L\'{e}vy bases.
As for the second part,
it is a simple consequence of Theorem~\ref{thm:isn_rm}, Lemma~\ref{lem:Miyoshi04a} and Proposition~\ref{prop:int_fld_meas}.
The third part follows from complete independence and Jensen's inequality.
\qed
\end{proof}
We shall now give some examples of $dcx$ ordered L\'{e}vy basis.
\begin{ex}
\label{ex:levy_basis} Let $\{x_i\}$ be a locally finite
deterministic configuration of points in $\mR^d$. Let $\{X^j_i\}_{i
\geq 1}, j=1,2$ be i.i.d sequence of infinite divisible random
variables such that $X^1_1 \leq_{cx} X^2_1$. (For example,  $X^1_1$
can be sum of two independent exponential r.v. with  mean $1/2$ and  $X^2_1$be
an exponential r.v. with mean $1$.) Define the L\'{e}vy bases as follows:
\[L_j(A) = \sum_{x_i \in A}X^j_i, \]
where $A$ is a bBs of $\mR^d$ and $j=1,2$. By Lemma~\ref{lem:lbasis}
and the fact that $X_1^1 \leq_{cx}X_1^2$
it follows that
$L_1 \leq_{dcx} L_2.$
\end{ex}
\begin{ex}\label{ex:Po-Levy}
Let $\tilde \Phi=\sum_{i}\varepsilon_{(x_i,Z_i)}$ be
an homogeneous Poisson p.p. on
$\mR^d$
independently marked by random variables $\{Z_i\}$
with mean $\lambda_0$.
Consider two random measures
$\Lambda_1 = \sum_{(x_i,Z_i)\in\tilde\Phi}\lambda_0\varepsilon_{x_i}$
and
$\Lambda_2 = \sum_{(x_i,Z_i)\in\tilde\Phi}Z_i\varepsilon_{x_i}$.
Note that $L_i$, $i=1,2$  are Levy basis. By
Lemma~\ref{lem:lbasis}  and the fact that $\lambda_0\leq_{cx} Z_i$,
conditioning on the number of points
and  using the same argument as in the proof of the second
statement  of Proposition~\ref{prop:mpp_idcx}
one can prove that $\Lam_1\leq_{dcx} \Lam_2$.
\end{ex}

\subsection{Poisson-Poisson Cluster Point Processes}
\label{sec:Po-Po-Cluster}
By Poisson-Poisson cluster p.p., we understand
a LCP with the Levy basis being a Poisson p.p.
This class deserves a separate mention due to the generality
of the ordering results that are possible. For rest of the section, assume that $h(x)$ is a non-negative measurable
function such that $\int_{\mR^d}h(x)dx = \lambda_0 < \infty.$

We shall now give an example of a parametric family of $dcx$
ordered Poisson-Poisson cluster  p.p..
Fix $\lam >0$. Let $\Phi_c, c > 0$ be a family of homogeneous
Poisson p.p.  on $\mR^d$ of intensity $c\lam$.
Let a non-negative function $h:\mR^d\times\mR^d\to\mR$ be given
and consider a family of shot noise fields
$\lambda_c(y)=\int_{\mR^d}\left(h(x,y)/c\right)\,\Phi_c(dx)$,
which are assumed a.s. locally Riemann integrable
with $\int_B\sE(\lambda_c(y))\,dy<\infty$ for bBs $B$.

\begin{prop}
\label{prop:par_lcp}
The family of shot-noise fields $\{\lambda_c(y)\}_{y\in\mR^d}$
is decreasing in $dcx$, i.e., for  $0<c_1\le c_2$ we have
$\{\lambda_{c_2}(y)\}\leq_{dcx}\{\lambda_{c_1}(y)\}$.
Consequently Cox($\lambda_{c_2}$)$\leq_{dcx}$\allowbreak
Cox($\lambda_{c_1}$).
\end{prop}

\begin{proof}
Note that $\{\lambda_c(x)\}$ can be seen as a shot-noise field
generated by the response function $h$ and the
Levy basis $L_c=(1/c)\Phi_{\lambda c}$.
By Lemma \ref{lem:lbasis} and  Theorem~\ref{thm:isn_rm},
it is enough to prove that $L_{c_2}(A) \leq_{cx} L_{c_1}(A)$ for $A$ bBs and
$c_2 > c_1 > 0$.

Since, $X \leq_{cx} Y$ implies that $aX \leq_{cx} aY$ for all
scalars $a > 0$, it suffices to prove that $L_{ca}(A) \leq_{cx} L_a(A)$ for
$A$ bBs and $c > 1, a > 0$. This essentially boils down to proving that
$N_{ca} \leq_{cx} c N_a, \ c > 1, a > 0$, where $N_{a}$ stands for a
Poisson r.v. with mean~$a$.
Let $\{X^n_i\}_{1 \leq i \leq n}$ and $\{Y^n_i\}_{1 \leq i \leq n},
\ n \geq 1$ be i.i.d. sequences of Bernoulli r.v's with
probability of success $ca/n$ and $a/n$, respectively, with $n\ge ca$.
Let $X^n = \sum_{i=1}^n X^n_i$ and $Y^n = \sum_{i=1}^n c Y^n_i$.
It is well known that $X^n, Y^n$ converge weakly to $N_{ca}, N_{a}$
respectively, as $n\to\infty$.
As convex order preserves weak convergence, we need to only prove that $X^n \leq_{cx} c Y^n$.
By the independence of summands,  it is enough to prove that
$X^n_i \leq_{cx} c Y^n_i$, which we  shall do in what follows.
Let $f$ be a convex and differentiable function.
Define $g(c) := \sE{f(X^n_i)}- \sE{f(cY^n_i)} = \frac{a}{n}\{c(f(1)-f(0))-f(c)+f(0)\}$.
Note that $g(1) = 0$. Hence, our proof is complete if we show that $g$ is decreasing in
$c > 1$. Indeed,
\begin{eqnarray*}
g^\prime(c) & = & \frac{a}{n}\{(f(1)-f(0)) - f^\prime(c)\} \\
& = & \frac{a}{n}\{f^\prime(b)-f^\prime(c)\} \leq 0, \ \ (b < c)
\end{eqnarray*}
where $b \in (0,1)$ by mean-value theorem and $f^\prime$ is increasing
due to convexity. \qed
\end{proof}

Poisson-Poisson cluster p.p. can be also $dcx$ compared to
a homogeneous Poisson p.p.. Let $\Phi$ and $\Phi'$ be homogeneous
Poisson p.p. with
intensities $\lambda < \infty$ and $\lambda \times \lambda_0$
respectively. Define $\mu(y) = \sum_{X_i \in \Phi} h(X_i-y)$.
Let $\Phi^{\prime \prime}$ be Cox($\mu(x)$).
\begin{prop}
\label{prop:po-po-cl}  Let $\Phi, \Phi'$, $\{\mu(y)\}$ be as above.
Assume that $\mu(y)$ is a.s. locally Riemann integrable
and $\sE(\mu(y))=\sE(\mu(0))<\infty$.
Then  $\Phi^\prime \leq_{dcx} \Phi^{\prime \prime}$.
\end{prop}

\begin{proof}
By the last statement
of Lemma~\ref{lem:lbasis} we have  $\lambda\,dx\leq_{dcx}\Phi(dx)$.
Note that $\lambda\times\lambda_0=\int_{\mR^d}h(x-y)\,
\lambda dx$ and thus by the second statement of
Theorem~\ref{thm:isn_rm} (note the assumption  $\sE(\mu(y))<\infty$)
$\{\lambda\times\lambda_0\}\leq_{dcx}\{\mu(y)\}$, where the $dcx$
smaller field is a deterministic, constant.
The result follows now from the second statement of
Lemma~\ref{lem:Miyoshi04a}
by assumption that $\mu(y)$ is a.s. Riemann integrable
and observing   that $\sE(\int_A\mu(y)\,dy)=\sE(\mu(0))\int_Ady<\infty$
for all bBs $A$.\qed
\end{proof}

\begin{rem}\label{rem:CoxPalm}
Consider Poisson p.p. $\Phi'$ and Cox($\mu$) as in
Proposition~\ref{prop:po-po-cl}. It is known that the Palm version
(given a point at the origin) of
$\Phi'$ can be constructed taking $\Phi'+\varepsilon_0$.
By~\cite[Proposition~2]{Moller03}, analogously,  Palm version of  Cox($\mu$)
can be taken as Cox$(\mu)+\varepsilon_0+\Phi''$,
where $\Phi''$ is an independent of
Cox$(\mu)$ Poisson p.p. with intensity $h(y-\xi)$ where
$\xi$ is sampled from the distribution $h(dx)/\int h(y)dy$.
This shows that one cannot expect $dcx$ ordering of the Palm versions
of $\Phi'$ and  Cox($\mu$).
\end{rem}

\subsection{Log Cox Point Processes}
\label{sec:log_cox}

This class of p.p. are defined by the logarithm of their
intensity fields.

An extension of LCP studied in \cite{Hellmund08} is Log-L\'{e}vy
driven Cox process (LLCPs). Under the notation of the previous
subsection, a p.p. $\Phi$ is said to be a LLCP if its
intensity field is given by
\[ \lambda(y) = \exp \left( \int_{\mR^d} k(x,y)L(dx) \right). \]
\cite{Hellmund08} allows for negative kernels and signed L\'{e}vy
measures but they do not fit into our framework. Suppose that $L_1
\leq_{idcx} L_2$, then $\Phi_1 \leq_{idcx} \Phi_2$ where $\Phi_i,i=1,2$
are the respective LLCPs of $L_i,i=1,2$ with kernel $k(.,.)$. These are
simple consequences of Theorem \ref{thm:isn_rm} and the exponential function being $icx$.

Another class is the Log-Gaussian Cox process (LGCPs)(see \cite{Moller98}).
A p.p. $\Phi$ is said to be a LGCP if its intensity field is
$\lambda(y) = \exp\{X(y)\}$ where $\{X(y)\}$ is a Gaussian random
field. Suppose $\{X_i(y)\}, i=1,2$ are two Gaussian random fields,
then $\{X_1(y)\} \leq_{idcx} \{X_2(y)\}$ if and only if  $\sE(X_1(y)) \leq
\sE(X_2(y))$ for all $y \in \mR^d$ and $\cov(X_1(y_1),X_1(y_2))
\leq \cov(X_2(y_1),X_2(y_2))$ for all $y_1,y_2 \in \mR^d.$ From the
composition rules of $idcx$ order, it is clear that $idcx$ ordering
of Gaussian random fields implies $idcx$ ordering of the
corresponding LGCPs. An example of
parametric $dcx$ ordered Gaussian random field is given
in~\cite[Sec~ 4]{Miyoshi04}.

\subsection{Generalized Shot Noise Cox Processes (GNSCPs)}
\label{sec:gnscp}

This class of Cox p.p. was first introduced and its various
statistics were studied in \cite{Moller05}. In simple terms, these
are Cox p.p. whose random intensity field is a shot-noise field
of a p.p. We say a Cox p.p. is GNSCP if the
random intensity field $\{\lambda(y)\}_{y \in \mR^d} $ driving the
Cox p.p. is of the following form : $\lambda(y) = \sum_j \gamma_j
k_{b_j}(c_j,y)$ where $(c_j,b_j,\gamma_j) \in \Phi$, a p.p. on
$\mR^d \times (0,\infty) \times (0,\infty)$. Also we impose the
following condition on the kernel $k$ : $k_{b_j}(c_j,y) =
\frac{k_1(c_j/b_j,y/b_j)}{b_j^d}$ where $k_1(c_j,.)$ is a density
with respect to the Lebesgue measure on $\mR^d$. We shall denote the GNSCP
driven by $\Phi$ as $\Phi^G$. This class includes
various known p.p. such as Neyman-Scott p.p., Thomas
p.p., Mat\'ern Cluster p.p. among others. The case when
$b_j$'s are constants and $\{(c_j,\gamma_j)\}$ is a Poisson p.p.
 is called as Shot Noise Cox process (See \cite{Moller03}).
Shot Noise Cox process are also LCPs. Suppose two p.p.
$\Phi_1 \leq_{dcx\resp{idcx;\,idcv}} \Phi_2$, then from Theorem
\ref{thm:isn_rm},
we infer that $\Phi^G_1 \leq_{dcx\resp{idcx;\,idcv}} \Phi^G_2$.

\subsection{Ginibre-Radii Like Point Process}
\label{sec:pp_rl}

Let $\{\Phi_i\}_{i \geq 0}$ be an i.i.d. family of p.p. on
$\mR^+$.
So, the points of each  p.p. $\Phi_i$
 can be sequenced based on their distance from the
origin. Let $\Phi$ be the p.p. formed {\em by picking the $i$th
point of $\Phi_i$}
for $i \geq 1$. We shall from now on abbreviate $\Phi([0,b])$ by $\Phi(b)$ for $b > 0$
and similarly for other p.p. used. Note the following representation for $\Phi(b)$
and $\Phi_0(b)$:
\[\Phi(b) =  \sum_{k \geq 1} \1[\Phi_k(b) \geq k] \ ; \ \Phi_0(b) =  \sum_{k \geq 1} \1[\Phi_0(b) \geq k]. \]
Let
\[\Phi^m(b) =  \sum_{k \geq 1}^m \1[\Phi_k(b) \geq k] \ ; \ \Phi^m_0(b) =  \sum_{k \geq 1}^m \1[\Phi_0(b) \geq k]. \]
By Lorentz's inequality (see \cite[Th.~3.9.8]{Muller02}), it follows that $(\Phi_1(b),\ldots,\Phi_m(b))
\leq_{sm} (\Phi_0(b),\ldots,\Phi_0(b))$,
where $sm$ stands for {\em supermodular} (see
\cite[\S~3.9]{Muller02}).
Define the $f:\mN^m \to \mR$ as follows : $f(n_1,\ldots,n_m) = \sum_{k \geq 1}\1[n_k \geq k]$.
It is easy to verify that both $f$ and $-f$ are $sm$ and
$f(n\wedge m) \leq f(n), f(m) \leq f(n \vee
m)$. In consequence  $g \circ f$ is $sm$ provided  $g$ is
$cx$ and
$ \sE(g(\Phi^m(b))) = \sE(g \circ f(\Phi_1(b),\ldots,\Phi_m(b))) \leq  \sE(g \circ f(\Phi_0(b),\ldots,\Phi_0(b))) =  \sE(g(\Phi^m_0(b)))$.
Hence $\Phi^m(b) \leq_{cx} \Phi_0^m(b)$ and using Lemma \ref{lem:weak_conv}, we get that $\Phi(b) \leq_{cx} \Phi_0(b)$.
To complete the proof $\Phi \leq_{dcx} \Phi_0$, one would require a multi-variate generalization of Lorentz's inequality which
we have been unable to prove.

We shall now explain the reasons for considering the above p.p.
$\Phi$. If we assume that $\Phi_i$ above  are Poisson,
then $\Phi$ is know to be a representation of the
p.p. of the squared radii
$|\Phi_G|^2=\{|X_n|^2 \ : \ X_n \in
\Phi_G\}$ of the {\em Ginibre process} $\Phi_G$
(see~\cite{Kostlan92,Ben06}).
It has been  observed  in simulations that this {\em determinental} p.p.
exhibits less clustering than the homogeneous Poisson p.p.
Our result can be seen as a first step
towards a formal statement of this property.

\remove{ For a point p.p. $\Phi$ on $\mR^d$,
denote by the point p.p. of radii on $\mR^+$. Let $\{\Phi_i\}_{i \geq 1}$ be an i.i.d. family of
homogeneous Poisson($1$) point p.p. on $\mR^+$ and $\Phi$ be as defined above. Let $\Phi^\prime, \Phi^{\prime \prime}$
be the Ginibre p.p. (\cite{Ben06}) and the homogeneous Poisson($1/\pi$) p.p. on $\mR^2$ respectively. Then it is
known that (\cite{Kostlan92}) $|\Phi^\prime|^2 = \Phi$ and $|\Phi^{\prime \prime}|^2 = \Phi_0$. Hence, the above result leads us to conjecture
that $|\Phi^\prime|^2 \leq_{dcx} |\Phi^{\prime \prime}|^2$ and more generally $\Phi^\prime \leq_{dcx} \Phi^{\prime \prime}$.}


\section{Applications to Wireless Communication Networks}
\label{sec:appln}

{}From the point of view of applications of our main result, what
remains is examples of interesting $dcx$ functions. In what follows,
we will provide such functions arising in the context of wireless networks. 
In many of the models we have assumed ordered point processes with i.i.d. marks. 
However due to Propnosition~\ref{prop:mark_cox}, the results hold for
independently  
marked Cox p.p. provided the respective intensity measures are ordered. 

\subsection{Coverage Process with Independent Grains}
\label{ssec:coverage}

The Boolean model $C(\Phi,G)$ defined earlier (see Definition~\ref{defn:boolean_model})
is the main object of analysis in the theory of Coverage processes (see~\cite{Hall88}).
The percolation properties of the Boolean model has been studied in~\cite{RMeester96}
while the connectivity properties of the Boolean model has been
studied in~\cite{Penrose03}.
For $\tilde \Phi$ as in the  Definition~\ref{defn:boolean_model}
of the Boolean model, denote by $V(y)=\sum_{(X_i,G_i)\in
  \tilde\Phi}\1[y \in X_i+ G_i]$ the number of grains
covering $y\in\mR^d$. Denote by $\psi(s_1,\ldots,s_n)$
the  {\em joint probability generating functional (p.g.f) of the number of
grains} covering locations $y_1,\ldots,y_n\in\mR^d$
$\psi(s_1,\ldots,s_n)=
\sE\Bigl(\prod_{j=1}^n s_j^{V(y_j)}\Bigr)$, $s_j\ge 0$, $j=1,\ldots,n$.
Note that  the function
$g(v_1,\ldots,v_n)=\prod_{j}^n s_j^{v_j}$ is $idcx$ when $s_j\ge1$ for all
$j=1,\ldots,n$ and is $ddcx$ when $0\le s_j\le 1$ for all $j$.

Thus the following result follows immediately from
Theorem~\ref{thm:isn_rm}, Proposition~\ref{prop:mpp_idcx} and Proposition \ref{prop:mark_cox}.
\begin{cor}\label{cor:coverage}
Let $\Phi_i$, $i=1,2$ be  a simple p.p. (of germs) on
$\mR^d$. Consider
the corresponding Boolean  models with the typical grain $G$
and, as above, denote the respective coverage number fields by $\{V_i(y)\}$ and
and their p.g.f  by $\psi_i$.
If $\Phi_1\leq_{dcx\resp{idcx;\,idcv}}\Phi_2$ then
$\{V_1(y)\}\leq_{dcx\resp{idcx;\,idcv}}\{V_2(y)\}$, with the result for $dcx$ holding
provided $\sE(V_i(y))<\infty$ for all $y$.
In particular, if $\Phi_1\leq_{idcx}\Phi_2$ then
$\sE(V_1(y)^\beta)\le\sE(V_2(y)^\beta)$ for all $\beta\ge1$.
If $\Phi_1\leq_{idcx\resp{ddcx}}\Phi_2$ then
  $\psi_1(s_1,\ldots,s_n)\le\psi_2(s_1,\ldots,s_n)$ for  $s_j\ge
  1$  (resp. $s_j\le   1$) $j=1,\ldots,n$.
\end{cor}
Note that $1-\psi(0,\ldots,0)$
represents the expected {\em coverage measure}, i.e.,the probability
whether the locations $y_1,\ldots,y_n$ are covered by at least
one grain.
In \cite[Section 3.8]{Hall88} it is shown that expected one-point coverage (or
volume fraction in case of stationary p.p.) for a stationary Cox p.p. and some
clustered p.p. is lower than that of a stationary, homogeneous
Poisson p.p..

Coverage processes arise in various applications.
In particular, in wireless communications
the points of the p.p. (germs)  usually represent locations of
antennas and their grains the respective communication regions.
In this context $V(y)$ is the number of antennas covering the point
$y$ and the coverage measure is the indicator that at least one of
them is able to reach $y$.
The application of the Boolean model to the modeling of wireless
communications dates back to the article of Gilbert~\cite{Gilbert61}
in~1961.

\subsection{Random Geometric Graphs (RGGs)}
This class of graphs has increasingly found applications in spatial networks.
For a detailed study of these graphs, see \cite{Penrose03}. A random geometric graph is defined as a graph with
$\Phi$ as the vertex set and the edge-set $E = \{\{X_i,X_j\}:|X_i-X_j| \leq r \}$. Clearly this is related to the Boolean model
defined in the previous subsection. One of the objects of interest in a RGG is the typical degree. Under the notation of the previous subsection,
the typical degree ($deg(\Phi,G)$) for a RGG formed by a stationary p.p. $\Phi$ and grain distribution $G$
is $deg(\Phi,G) = \frac{1}{\lam |A|}\sum_{X_i,X_j \in \Phi} \1[X_i \in A] \1[X_i \neq X_j] \1[(X_i+G_i)\cap(X_j+G_j) \neq \emptyset]$,
where $A$ is a bBs. If $G = B_0(r), r > 0$, then $\sE(deg(\Phi,G)) = K(r)$ is the {\em Ripley's $K$~function} defined in Section \ref{sec:moments}.
The following result follows easily from Theorem \ref{thm:isn_rm}, Proposition \ref{prop:product_space} and Proposition \ref{prop:mark_cox}.

\begin{cor} \label{cor:rgg}
Suppose that simple p.p. $\Phi_1 \leq_{dcx} \Phi_2$, then  $deg(\Phi_1,G) \leq_{idcx} deg(\Phi_2,G)$.
\end{cor}

\subsection{Interference in Wireless Communications}
\label{ssec:conn_sinr}

The Boolean model is not sufficient for analyzing wireless networks
as it ignores the fact that in  radio communications
signal received from one particular transmitter  is jammed by the
signals received from the other transmitters.
According to information theory as well as existing technology,
the quality of a given radio communication link
is determined by the so called  {\em signal to interference and noise
ratio} (SINR) at the receiver of this link.
{} a  mathematical point of view,
the interference in the above considerations
is just the sum of the powers of the signals received from all
transmitters (perhaps except own transmitter(s)).
It is then the shot-noise field of received powers
that  plays important role in determining the
connectivity and the capacity of the network in a broad sense.
The foundations of the theory of SINR coverage processes
are quite recent (see \cite{GK00,Baccelli01,Dousse_etal2006,Baccelli06}).
In what follows, we shall study the impact of
structure of the p.p. of interferers on given radio links.

 Consider a set of $n$
emitters $\{x_i\}$ and $n$ receivers $\{y_i\}$.
Suppose that the signal received by $y_i$ from  $x_k$ is $S_{ki}$. These
$\{S_{ik}\}$ are assumed to be independent.
The assumption of independence is
due to the phenomenon of {\em fading}.
Let the set of additional
interferers be modeled  by a i.i.d. marked p.p. $\tilde\Phi =
\varepsilon_{(X_j,(Z_j^1,\ldots,Z_j^n)}$, independent of $\{S_{ik}\}$,
where $Z_j^n$ is the power received by the receiver $y_i$ from the
interferer located at $X_j$.
Denote the {\em background noise} random variable by
$W$.

We say that the  signal from $x_i$ is successfully
received by $y_i$ if
$S_{ii}/(W+I_i+V_i)>T $
where $I_i=\sum_{k\not=i}S_{ki}$ and
$V_i=\sum_{j}Z_j^i$ is the interference received at $y_i$ from
the set of other emitters $\{x_k:k\not=i\}$
and interferers in $\tilde\Phi$, respectively, and  $T>0$ is some
(assume constant)
required  {\em SINR  threshold}.
If we denote by $p$, the probability
of successful reception  of signals from each $x_i$ to $y_i$, then
\begin{eqnarray}
p & = & \sP(S_{ii} > (W+I_i+V_i)T\quad  \forall i=1,\ldots,n) \no \\
& = & \sE(\prod_i \overline  F_{ii}(T(W+I_i+V_i)))\,, \label{eqn:ddcx_fn}
\end{eqnarray}
where $\overline F_{ii}(s)=\sP(S_{ii}\ge s)$ and the second equality
is due to independence.
Given  $\{I_i:i=1,\ldots,n\}$ and $W$,
the expression under expectation in~(\ref{eqn:ddcx_fn})
can be viewed as a function
of the value of the shot-noise vector  $(V_1,\ldots,V_n)$
evaluated with respect to $\tilde\Phi$.
Theorem \ref{thm:isn_rm} and Proposition~\ref{prop:mark_cox}
implies the following result concerning the impact of the structure of the
 set of interferers on $p$.
\begin{cor}\label{cor:coverage-sinr}
Consider emitters $\{x_i\}$, receivers $\{y\}_i$, powers $\{S_{ki}\}$
as above. Let $\tilde\Phi_u$, $u=1,2$ be two simple marked p.p. of interferers.
Denote by $p_u$, $u=1,2$ the probability of successful reception
given by~(\ref{eqn:ddcx_fn}) in the model with the set of interferers
$\tilde\Phi_u$.
Assume the product of tail distribution functions of the received powers
$\prod_{i=1}^n\overline F_{ii}(s_i)$ be $dcx$.
If $\Phi_1\leq_{ddcx}\Phi_2$ then $p_1\le p_2$.
\end{cor}

It is quite natural to assume $ddcx$ $\prod_{i=1}^n\overline F_{ii}(s_i)$.
For example the constant emitted power~$P$, omni-directional path-loss
function $l(r)$ and  {\em Rayleigh fading} in the radio channel
implies $S_{ki}=PH_{ki}/l(|x_k-y_i|)$, where $|H_{ki}|$ are i.i.d.
exponential random variables with mean~$1$. In this case
$\prod_{i=1}^n\overline F_{ii}(s_i)$ is
ddcx.~(${}^*$)\footnote{(${}^*$)
Recently in~\cite{Ganti08}, under the assumption of Rayleigh fading,
direct analytical methods have been used
to compare the probability of successful
reception in  Poisson p.p. and a class of Poisson-Poisson cluster p.p.
known as Neyman-Scott p.p. for both stationary and Palm
versions. These results relay on explicit expressions for this
probability known in the considered cases. Further, it is shown that for a certain choice of parameters,
Palm version of the Poisson-Poisson cluster p.p. has a worser 
probability of successful reception than the Poisson p.p.. 
In our terminology, it simply means that 
the corresponding Palm versions aren't $ddcx$ ordered as the
connectivity probability is a $ddcx$  
function (Eqn.~\ref{eqn:ddcx_fn}) of the integral shot-noise fields of
the corresponding Palm versions. This strengthens
Remark~\ref{rem:CoxPalm} by showing that $idcx$ ordering of
Palm versions is the best one can obtain in full generality. }


\section{Conclusions and Open Questions}
\label{sec:concl}

To the best of our knowledge,
this is the first study of $dcx$ ordering of random measures and p.p..
We have defined the $dcx$ order and characterized it  by
finite dimensional distributions of the measure values on disjoint bBs
of the space. As the main result, we have proved that
the integrals of some non-negative kernels
with respect to  $dcx$ ordered random measures inherit this ordering
from the measures.  This was shown to be a very useful tool in study
of many particular characteristics of random measures and in the
construction and analysis of stochastic models.

In this paper, we have also left  several open questions. Here we
briefly summarize them.

\begin{itemize}
\item Our $dcx$ order is defined via finite dimensional distributions
  of random measures. This makes the verification of $dcx$ order more
  easy but requires additional work when studying functionals,  which
  cannot be explicitly  expressed in terms of the values of the measure
  on some finite collection  of bBs as, e.g., an  integral of the measure.
Considering a  $dcx^{1*}$ order on the space of measures could facilitate
the former task.  However, the precise regularity conditions of the
$dcx^{1*}$ functional on the space of measures which would guarantee
the equivalence between these two approaches are not known
(cf Section~\ref{sec:alt_defn}).
\item Comparisons of Ripley's functions (see Proposition~\ref{prop:Ripley})
and pair correlation functions (Corollary~\ref{cor:pcf})
seem to indicate that the higher in  $dcx$ order processes cluster
their points more. We have shown examples of p.p.,
which are larger than Poisson one, namely Cox p.p., which
indeed exhibit more clustering than in Poisson p.p..
It would be interesting to show  examples of p.p. which are
$dcx$ smaller  than Poisson one, and which exhibit
less clustering than it. Mat\'ern ``hard core'' p.p. and  Ginibre p.p.
are some natural candidates for this.
\item We have studied  $dcx$ order
that takes into account the dependence structure
and the variability of the marginals or random measures.
It seems plausible to study in a similar manner
other orders such as convex, component-wise convex order
etc. Note however that the supermodular order does not seem to be
a reasonable one in the context of random measures. The reason is that
it allows to compare only measures with the same finite dimensional
distributions, and thus a Poisson p.p. can only be (trivially) compared
in this order to itself. Indeed, Poisson finite dimensional
distributions imply total independence property and thus uniquely
characterize Poisson p.p. (cf~\cite[Lemma~2.3.I]{DVJ88}).


\end{itemize}
%

\section{Appendix}
\label{sec:Appendix}



In order to make the paper more self-contained,
we shall recall now  some basic results on stochastic
orders used in the main stream of the paper.
The following two lemmas can be found
in \cite[Chapter 3]{Muller02}.

\begin{lem}
\label{lem:dcx_diff}
\begin{enumerate}
 \item A twice differentiable function $f$ is directionally convex if and only if
$\dfrac{\partial^2}{\partial x_i \partial x_j}f(x) \geq 0,$
for all $x, 1 \leq i,j \leq n.$
\item The stochastic order relation $\leq_{dcx}$ is generated by infinitely differentiable
$dcx$ functions.
\end{enumerate}
\end{lem}

Due to the above lemma, at some places we only prove that two random
vectors are ordered with respect to
twice differentiable $dcx$ functions and conclude that they are $dcx$
ordered.

We denote by  $\cvgdist$ convergence in distribution (weak
convergence).
\begin{lem}
\label{lem:weak_conv}Let $(X^{(k)}:k=1,\ldots)$ and $(Y^{(k)}:k=1,\ldots)$
be sequences of random vectors. Suppose $X^{(k)} \leq_{dcx} Y^{(k)}$ for all $k \in \mN$. If $X^{(k)} \cvgdist
X$ and $Y^{(k)} \cvgdist Y$  and if moreover $\sE(X^{(k)}) \rar \sE(X)$ and $\sE(Y^{(k)}) \rar \sE(Y)$, then
$X \leq_{dcx} Y$.
\end{lem}
The following result is from Lemmas 2.17 and 2.18 of
\cite{LMeester93}.
\begin{lem}
\label{lem:LMeester93}
\begin{enumerate}
\item For $i=1,\ldots,m$ let $(S_j^i:j=1,\ldots)$ be independent
  sequences of i.i.d. non-negative random variables.
Suppose $f$ is $dcx\resp{idcx;\,idcv}$, then
$g(n_1,\ldots,n_m) =
\sE(f(\sum_{j=1}^{n_1}S^1_j,\ldots,\sum_{j=1}^{n_m}S^m_j))$ is also
$dcx\resp{idcx;\,idcv}$.

\item Let $N_i, i = 1,\ldots,k$ denote $k$ mutually independent Poisson
r.v. where the mean of $N_i$ is $\lambda_i$. If $\phi:
\mN^k \rar \mR$ is $dcx\resp{idcx;\,idcv}$, then
$g(\lambda_1,\ldots,\lambda_k) = \sE(\phi(N_1,\ldots,N_k))$ is also
$dcx\resp{idcx;\,idcv}$.
\end{enumerate}
\end{lem}
The first part of the following lemma is an easy extension of the
one-dimensional version in~\cite{LMeester93}.
The second part, which we prove in what follows, is a further extension of it.

\begin{lem}
\label{lem:Miyoshi04a} Suppose $\{X(s)\}_{s \in \mR^d}$ and
$\{Y(s)\}_{s \in \mR^d}$ are two non-negative real-valued and a.s.
locally Riemann
integrable  random fields.
For some $n\ge1$ and  disjoint bBs  $I_1,\ldots,I_n$
denote $J_X^i=\int_{I_i}X(s)ds$, $J_Y^i=\int_{I_i}Y(s)ds$.
\begin{enumerate}
\item If $\{X(s)\} \leq_{idcx\resp{idcv}} \{Y(s)\}$, then
$(J_X^1,\ldots,J_X^n)\leq_{idcx\resp{idcv}}(J_Y^1,\ldots,J_Y^n)$.
%
%
for any $n$ and for any $I_1,\ldots,I_n$ disjoint bBs.

\item
Suppose further that $\sE(\int_{A}X(x)dx) < \infty$ for all bBs $A$ in $\mR^d$
and similarly for $\{Y(x)\}$. If $\{X(x)\}\leq_{dcx}(dcv)
\{Y(x)\}$,
then $(J_X^1,\ldots,J_X^n)\leq_{dcx(dcv)}(J_Y^1,\ldots,J_Y^n)$.

\remove{\item Suppose $\{X(s)\} \leq_{dcx(idcx)}\{Y(s)\}$.
For any $f$ $idcx(idcx^+)$
\[
\sE\left(X(0)f(J_X^1,\ldots,J_X^n)\right)
\leq \sE\left(Y(0)f(J_Y^1,\ldots,J_Y^n)\right)\,.\]
}
\end{enumerate}

\end{lem}

\remove{\begin{rem}
The converse of the first part of the lemma can be proved under some
regularity conditions on the random fields. What we need precisely
is the following: 1)~$X(s) = \lim_{\epsilon \to
0}\frac{1}{\|B(s,\epsilon)\|}\int_{B(s,\epsilon)}X(s)ds$ in
distribution where $B(s,\epsilon)$ is the ball of radius $\epsilon$
centered at $s$. 2)~The convergence holds also in  expectations.
Then the result shall follow as a consequence of Lemma~\ref{lem:weak_conv}.
For stationary random fields, condition 2 is
trivially true and the equality in condition 1 is true a.s. under
some regularity assumptions (for ex. a.s. continuity of sample paths) on the random fields.
We shall see some examples in the next section.
\end{rem}

\begin{rem}
The third part of the lemma was stated
in \cite{Miyoshi04} (See Lemma 5) for all $idcx$ functions.
The erroneous proof used the fact that
a product of $idcx$ functions is $idcx$. This is true only for
positive $idcx$ functions.
\end{rem}
}
\begin{proof}
$(2)$  We shall prove for $d=1$ and as can be seen from the
proof, the generalization is fairly straightforward.

We need to prove that $(\int_{I_1}X(s)ds,\ldots,\int_{I_n}X(s)ds) \leq_{dcx}
(\int_{I_1}Y(s)ds,\ldots,\int_{I_n}Y(s)ds),$ for
$I_i,i=1,\ldots,n$ disjoint bBs. We shall give an
approximation satisfying the assumptions of Lemma
\ref{lem:weak_conv}. Let $I_i = [a_i,b_i]; a_i,b_i \in \mR
,i=1,\ldots,n$. Let $\{(t_{mj}^i)_{1 \leq j \leq k_m},
i=1,\ldots,n\}$ be the sequences of $m$th nested partition of each
interval. The middle Riemann sum can be given as follows :
$X^m(I_i) = \sum_j X(t_{mj}^i)(t_{m(j+1)}^i-t_{mj}^i), i=1,\ldots,n , k \in
\mN$ and similarly for $Y(x)$. These are the variables satisfying the approximation as
in Lemma \ref{lem:weak_conv}. As $X(s)$ is Riemann integrable,
\[ (X^m(I_1),\ldots,X^m(I_n)) \rar (J_X^1,\ldots,J_X^n) \]
a.s. and hence in distribution. It is also clear the middle Riemann sums of $X(\cdot)$ and
$Y(\cdot)$ are ordered. What remains to prove is that $\sE{X^m(I_i)} \rar \sE{J_X^i}$.
In the last term, by Fubini, we can interchange the expectation and integral and hence it
suffices to prove $\sE{X^m(I_i)} \rar \int_{I_i}\sE{X(s)}ds$. Our assumption implies that
this is true.
\remove{
$(3)$ Observe the following : If $f: \mR^n \to \mR$ is $idcx(idcx^+)$ then $g(x_0,x)= x_0^+f(x):\mR^{n+1} \to \mR$
is $dcx(idcx)$, where $a^+$ for $a \in \mR$ denotes the non-negative part of $a$. Check that if $f$ is twice differentiable,
then so is $g$ for $x_0 \geq 0$ and all the second partial derivatives of $g$ are non-negative. Thus by first part of Lemma \ref{lem:dcx_diff} ,
$dcx$ condition is satisfied for $x_0 \geq 0$. For $x_0 < 0$, the $dcx$ condition can be checked directly. Second part of Lemma \ref{lem:dcx_diff}
completes the proof.} \qed
\end{proof}

\ack DY was supported by a grant from EADS, France. DY also wishes to thank Prof. Fred Huffer for
sharing the technical report.

\renewcommand{\AAP}{\emph{Adv. Appl. Probab.}}
\renewcommand{\JAP}{\emph{J. Appl. Probab.}}
\renewcommand{\AP}{\emph{Ann. Probab.}}

\end{document}